\documentclass[11pt, reqno]{amsart}

\DeclareMathAccent{\mathring}{\mathalpha}{operators}{"17}

 \usepackage{color}

\newcommand{\mysection}[1]{\section{#1}
      \setcounter{equation}{0}}

\newtheorem{theorem}{Theorem}[section]
\newtheorem{lemma}[theorem]{Lemma}

\theoremstyle{definition}
\newtheorem{assumption}{Assumption}[section]

\theoremstyle{remark}
\newtheorem{remark}{Remark}[section]

\newcommand\cbrk{\text{$]$\kern-.15em$]$}}
\newcommand\opar{\text{\raise.2ex\hbox{${\scriptstyle | }$}\kern-.34em$($} }

 \makeatletter
 \def\dashint{%
 \operatorname%
 {\,\,\text{\bf--}\kern-.98em\DOTSI\intop\ilimits@\!\!}}
 \makeatother

 \newcommand{\WO}{\overset{ 
 \scriptscriptstyle0}%
 { W}\,\!}

\newcommand\bR{\mathbb{R}}

\newcommand\bZ{\mathbb{Z}}

\newcommand\cF{\mathcal{F}}

\newcommand\cP{\mathcal{P}}

\newcommand\cH{\mathcal{H}}

\newcommand\cS{\mathcal{S}}

\newcommand\dist{{\rm dist}\,}

\begin{document}

\title[Existence for fully nonlinear elliptic
equations]
{On the existence of smooth 
solutions for fully nonlinear elliptic
equations with measurable ``coefficients''
without convexity assumptions}

\author{N.V. Krylov}
\thanks{The  author was partially supported by
  NSF Grant DMS-1160569}
\email{krylov@math.umn.edu}
\address{127 Vincent Hall, University of Minnesota,
 Minneapolis, MN, 55455}



\begin{abstract}
We show that for any uniformly elliptic
fully nonlinear second-order equation with bounded 
measurable ``coefficients''
and bounded ``free'' term
one can find an approximating
equation which has a unique
 continuous and having the second derivatives
locally bounded  solution  
in a given smooth domain  with  
smooth boundary data. The approximating
equation is constructed in such a way
that it modifies the original one
only for large values of the unknown
function and its derivatives.
\end{abstract}

\maketitle

\mysection{Introduction and  main result}

                                             \label{section 2.5.1}

In this article, we consider elliptic  equations
\begin{equation}
                                                \label{7.29.1}
H[v](x):= H(v( x),D v( x),D^{2}v( x), x)=0
\end{equation}
in subdomains of $\bR^{d } $, where
$$
\bR^{d}=\{x=(x_{1},...,x_{d}):x_{1},...,x_{d}\in
\bR=(-\infty,\infty)\}.
$$
Here
$$
D^{2}u=(D_{ij}u),\quad Du=(D_{i}u),\quad
D_{i}=\frac{\partial}{\partial x_{i}},\quad D_{ij}=D_{i}D_{j}.
$$
We introduce $\cS$ as the set of symmetric $d\times d$ matrices,
fix a constant 
$\delta\in(0,1]$, and  set
$$
\cS_{\delta}=\{a\in\cS:\delta|\xi|^{2}\leq a_{ij}\xi_{i}\xi_{j}
\leq\delta^{-1}|\xi|^{2},\quad \forall \xi\in\bR^{d}\},
$$
where and everywhere in the article the summation convention is enforced
unless specifically stated otherwise.

Recall that   Lipschitz
continuous functions are almost everywhere differentiable.

\begin{assumption}
                                    \label{assumption 9.23.1}
(i) The function $H(u,x)$,
$u=(u',u'')$, 
$$
u'=(u'_{0},u'_{1},...,u'_{d})
\in\bR^{d+1},\quad u''\in\cS,\quad x\in\bR^{d},
$$
 is measurable with respect to $x$ for any $u$,
and Lipschitz continuous
in $u$ for every $x\in\bR^{d}$. 

(ii) For any $x$, at all 
points of differentiability of $H(u,x)$ with respect to $u$,
we have
$$
(H_{u''_{ij}})\in \cS_{\delta},\quad
|H_{u'_{k}}|\leq \delta^{-1},\quad k=1,...,d,
\quad 0\leq-H_{u'_{0}}\leq  \delta^{-1}.
$$

(iii)
Finally,  
$$
\bar{H}:= \sup_{x
\in\bR^{d}}|H(0,x)|<\infty. 
$$
 
\end{assumption}

Let $\Omega$ be an open bounded subset of $\bR^{d}$
with $C^{2}$ boundary and take a function $g\in C^{1,1}
(\bar{\Omega})$.
Here is our   main result, in which $K\geq0$
is a fixed constant.

\begin{theorem}
                                    \label{theorem 9.23.1}
There is a  constant  $\hat{\delta}\in(0,\delta]$
depending only on $\delta$ and $d$ and there exists
a function $P(u) $ (independent of $x$), 
satisfying Assumption
\ref{assumption 9.23.1} 
with $\hat{\delta}$   
in  place of 
$\delta$, 
such that   the equation
\begin{equation}
                                               \label{9.23.2}
\max(H[v],P[v]-K)=0
\end{equation}
in $\Omega$  (a.e.) with boundary condition $v=g$ on $\partial\Omega$
has a unique solution 
$v\in C^{0,1}(\bar{\Omega})\cap C^{1,1}_{loc}(\Omega)$.
In addition,  for all $i,j$, and $p\in(d,\infty)$,
\begin{equation}
                                                \label{1.13.1}
|v|,|D_{i}v|,\rho|D_{ij} v |\leq N(\bar{H}+K
+\|g\|_{C^{1,1}(\Omega)})\quad\text{in}
\quad \Omega \quad (a.e.), 
\end{equation}
\begin{equation}
                                                \label{1.13.2}
\|v\|_{W^{2}_{p}(\Omega)}\leq
N_{p}(\bar{H}+K+\|g\|_{W^{2}_{p}(\Omega)}),
\end{equation}
\begin{equation}
                                                \label{2.28.1}
\|v\|_{C^{\alpha}(\Omega)}\leq N
(\|H[0]\|_{L_{d}(\Omega)}+\|g \|_{C^{\alpha}(\Omega)}),
\end{equation}
where 
$$
\rho(x)=\dist(x,\bR^{d}\setminus\Omega),
$$
$\alpha\in(0,1)$ is a constant depending
only on $d$ and $\delta$,
   $N$ is a constant depending only on $\Omega$ 
 and $\delta$, whereas $N_{p}$ only depends on the same objects and
$p$.

Finally, $P(u )$ is
constructed on the sole basis of $\delta$ 
and $d$, it is  positive homogeneous of degree one
and convex in $u$.

\end{theorem}
\begin{remark}
                                   \label{remark 1.13.1}
If we drop \eqref{1.13.2} and replace $C^{1,1}$
in \eqref{1.13.1} with $C^{1,\alpha}$, $\alpha\in(0,1]$,
 the assumptions
of Theorem \ref{theorem 9.23.1} about smoothness 
of $\Omega$ and $g$
can be somewhat relaxed. It is sufficient to
have the exterior ball condition on $\Omega$
and $g\in C^{1,\alpha}(\Omega)$.
Furthermore, if we multiply the derivatives
in \eqref{1.13.1} by one more $\rho$, then one can deal with $\Omega$
such that, for each boundary point $x_{0}$ and all
 $r>0$ small enough with the smallness
 independent of $x_{0}$, there is a ball of radius
$\varepsilon r$ at the distance $r$ from $x_{0}$
lying outside $\Omega$. Here $\varepsilon>0$ is a fixed constant.
Of course, in that case the asserted regularity  
should be changed to $v\in C^{0,\beta} (\bar{\Omega})
\cap C^{1,1}_{loc}(\Omega)$, where $\beta\in(0,1]$
is determined by other parameters of the problem.
 All these and other possible extensions 
and generalizations are left
to the interested reader.
\end{remark}

To the best of the author's knowledge Theorem \ref{theorem 9.23.1} is
the first uniqueness and {\em existence\/} result for {\em general\/}
fully nonlinear elliptic equations with
{\em measurable\/} coefficients without {\em convexity\/}
assumptions. In case $H$ is 
Lipschitz continuous
in $x$ the theory of viscosity solutions
provides the existence and uniqueness.
 Generally, one only knows that
such solutions are in $C^{1+\alpha}$ (see Trudinger
\cite{Tr89}). 
N. Nadirashvili and S. Vl\v{a}dut \cite{NV}
  found an example in which
viscosity solutions even for $H$ independent of $x$
do not have bounded second-order derivatives.

It is also worth  
mentioning that
M. G. Crandall, M. Kocan, and A. \'Swi{\c e}ch 
\cite{CKS00} developed a theory of $L_{p}$-viscosity
solutions for equations with measurable coefficients
(see also the references therein).

As far as a priori estimates in Sobolev spaces are concerned,
L. Caffarelli was the first author who derived 
interior $W^{2}_{p}$ estimates
under an assumption that certain estimates
hold for equations with zero ``free'' term,
  which are known to hold only for $H$ that are
either convex or concave with respect to
$v,Dv,D^{2}v$
(see   \cite{Caf89} and
\cite{CC95}). A particular case of $C^{2+\alpha}$ a priori
estimates without this assumption is presented in \cite{CC03}.
Another case is found in \cite{Ko09}.

  The activity which started in \cite{Caf89}
was continued by L. Wang in \cite{Wa92}
who obtained similar interior
a priori estimates for parabolic equations,
by M. G. Crandall, M. Kocan, and A. \'Swi{\c e}ch 
\cite{CKS00} who established the {\em solvability\/} 
in local Sobolev spaces of
the boundary-value problems
for fully nonlinear parabolic equations, and by  N.~Winter
\cite{Wi09} who established the solvability in the global
$W^2_p$-space of the associated boundary-value problem
in the elliptic case. In the existence parts in 
\cite{CKS00} and \cite{Wi09} the function
 $H$ is supposed to be convex with respect to
$ D^{2}v$ and continuous in $x$
(concerning the latter assumption see
\cite[Remark 2.3]{Wi09}, \cite{Kr10}, and
 \cite[Example 8.3]{CKS00}).
  However, in 
the above references the authors consider equations like
\eqref{7.29.1}  with the right-hand side which is not zero
but rather a function from an $L_{p}$-space. In our setting
we can only treat bounded right-hand sides.
 
 Recently a new method, very different
from the methods in the above cited references,
 emerged in \cite{DKL}
for treating fully nonlinear
elliptic and parabolic equations with VMO ``coefficients''.
Still the convexity of $H$ with respect to $D^{2}v$
is required in \cite{DKL} while proving the existence result.

In our
Theorem \ref{theorem 9.23.1} we do not impose any
convexity assumption on $H$ and allow it to be just
measurable in $x$. By the way, this theorem
 is obviously applicable to {\em
linear\/} equations. Yet we  approximate them with nonlinear
ones. 

The methods of the present article are
quite elementary and, apart from what is related to
  \eqref{1.13.2} and \eqref{2.28.1} and uniqueness,
 are not  using anything
from any existing theory of partial differential
equations in the main case where
$H$ depends only on pure second-order derivatives
and is continuous in $x$. Our main tool is finite-difference
approximations,  best demonstrated
in Sections \ref{section 9.22.1} 
and \ref{section 10.18.2}, which the reader may like
to read first.

\begin{remark}
                                    \label{remark 12.21.1}
It is almost obvious that
Assumption \ref{assumption 9.23.1} (ii) 
is equivalent to the requirement
that, for any $u\in\bR^{d+1}\times\cS$, $x,\xi\in\bR^{d}$,
$\eta\in\{\pm e_{1},...,\pm e_{d}\}$, where $e_{1},...,e_{d}$
 is the set
of standard basis vectors in $\bR^{d}$, and $r\geq0$, we have
$$
\delta|\xi|^{2}\leq H(u',u''+\xi\xi^{*},x)-
H(u',u'' ,x)\leq \delta^{-1}|\xi|^{2},
$$
$$
|H(u'+r(0,\eta),u'' ,x)-H(u',u'' ,x)|\leq \delta^{-1}r,
$$
$$
 H(u',u'' ,x)-\delta^{-1}r
\leq H(u'+r(1,0),u'' ,x)\leq H(u',u'' ,x) ,
$$
where $(0,\eta)=(0,\eta_{1},...,\eta_{d})$ and
$(1,0)=(1,0,...,0)$.

\end{remark}

\begin{remark}
                                    \label{remark 2.28.1}
Estimate \eqref{2.28.1} follows from
other assertions of Theorem \ref{theorem 9.23.1}
and the classical results about linear equations
with measurable coefficients (see, for instance,
Section 9.9 of \cite{GT}). Indeed, as is easy to see
for $v\in W^{2}_{p}(\Omega)$ satisfying
\eqref{9.23.2} we have that
$$
-\max(H[0],P[0]-K)=\max(H[v],P[v]-K)-\max(H[0],P[0]-K)
$$
$$
=a_{ij}D_{ij}v+b_{i}D_{i}v-cv
$$
with some functions $a=(a_{ij})\in\cS_{\hat{\delta}}$,
$|b_{i}|\leq\hat{\delta}^{-1}$, $0\leq
c\leq\hat{\delta}^{-1}$  (cf. the proof of Lemma
\ref{lemma 10.4.2}). Furthermore,
$|\max(H[0],P[0]-K)|\leq|H[0]|$.

 The assertion of Theorem \ref{theorem 9.23.1}
concerning uniqueness in our class of functions
is also a classical result derived from
the Alexandrov estimate.
\end{remark}

Here is
an almost trivial generalization of Theorem
\ref{theorem 9.23.1} which may be useful in some
applications.
\begin{theorem}
                                  \label{theorem 2.20.1}
Let $\phi\in C^{2}(\bar{\Omega})$ be a 
strictly positive function. Then
there is a  constant  $\hat{\delta}\in(0,\delta]$
depending only on $\delta$, $\phi$, and $d$ and there exists
a function $P(u) $ (independent of $x$), 
satisfying Assumption
\ref{assumption 9.23.1} 
with $\hat{\delta}$   
in  place of 
$\delta$, 
such that all assertions of Theorem  \ref{theorem 9.23.1}
hold true if we replace $P[u]$ with $P[\phi u]$
and allow the constants to also depend on $\phi$.
\end{theorem}

This result is obtained from Theorem
\ref{theorem 9.23.1} just by replacing there $u$ and $g$
with $\phi u$ and $\phi g$, respectively.

Here is a version of Theorem \ref{theorem 9.23.1}
which is obtained by just replacing $H(u,x)$
with $-H(-u,x)$.

\begin{theorem}
                                    
With $P$ from Theorem \ref{theorem 9.23.1}
the equation
$$
\min(H[v],-P[-v]+K)=0
$$
in $\Omega$  (a.e.) with boundary condition $v=g$ on $\partial\Omega$
has a unique solution 
$v\in C^{0,1}(\bar{\Omega})\cap C^{1,1}_{loc}(\Omega)$.
In addition,  for all $i,j$, and $p\in(d,\infty)$,
$$
|v|,|D_{i}v|,\rho|D_{ij} v |\leq N(\bar{H}+K
+\|g\|_{C^{1,1}(\Omega)})\quad\text{in}
\quad \Omega \quad (a.e.), 
$$
$$
\|v\|_{W^{2}_{p}(\Omega)}\leq
N_{p}(\bar{H}+K+\|g\|_{W^{2}_{p}(\Omega)}),
$$
$$
\|v\|_{C^{\alpha}(\Omega)}\leq N
(\|H[0]\|_{L_{d}(\Omega)}+\|g \|_{C^{\alpha}(\Omega)}).
$$
where 
$\alpha $,   $N$, and  $N_{p}$ are the constants
from Theorem \ref{theorem 9.23.1}.

\end{theorem}

It is an interesting issue as to what is happening to
$v=v_{K}$ as $K\to\infty$, where $v_{K}$ is taken from Theorem 
\ref{theorem 9.23.1}. We have the following

{\em Conjecture\/}. Assume that $H(u,x)$ is Lipschitz
continuous with respect to $x$ with Lipschitz constant
equal to a constant times $1+|u|$. Let $w$ be a unique viscosity
solution of equation \eqref{7.29.1} in $\Omega\in C^{2}$
with boundary condition $g\in C^{3}$. Then $|w-v_{K}|\leq N/K$
where $N$ is a constant.

To conclude our comments about Theorem \ref{theorem 9.23.1}
 we show how $P$
is constructed.
By Theorems 3.1 of \cite{Kr11} there exists
a set 
$$
\{ l_{1},...,l_{m}\}
\subset \bZ^{d },
$$ 
$m=m(\delta,d)\geq d$, chosen on the sole
 basis of
knowing $\delta$ and $d$ and 
 there exist
a constant  
$$
\hat{\delta}=\hat{\delta}(\delta,d ) 
\in(0,\delta/4]
$$
such that:

(i) We have 
$$
e_{i},e_{i}\pm e_{j}\in \{l_{1},...,l_{m}\}
=\{-l_{1},...,-l_{m}\}
$$
for all $i,j=1,...,d$  (recall that $e_{1},...,e_{d}$ is the 
 standard orthonormal basis of $\bR^{d}$);

(ii) There   exist
real-analytic functions $\lambda_{1}(a),...,
\lambda_{m}(a)$ on $\cS_{\delta/4}$ such that
for any $a\in\cS_{\delta/4}$
\begin{equation}
                                         \label{4.8.01}
a\equiv\sum_{k=1}^{m}\lambda_{k}(a)l_{k}l_{k}^{*},
\quad \hat{\delta}^{-1}\geq\lambda_{k}(a)\geq\hat{\delta}
,\quad \forall k.
\end{equation}

 Now introduce
$$
\cP(z)  =\max_{\substack{\hat{\delta}/2\leq a_{k}\leq
2\hat{\delta}^{-1} \\k=1,...,m} }
\max_{\substack{ |b_{k}|\leq
2\hat{\delta}^{-1}\\k=1,...,d_{1}} }
\max_{\hat{\delta}/2 \leq c\leq
2\hat{\delta}^{-1}}\big[\sum_{k=1}^{m} a_{k}
z''_{k}+\sum_{k=1}^{d}b_{k}z'_{k} -cz'_{0}\big],
$$
and for $u=(u',u'')\in\bR^{d+1}\times\cS$ define
$$
P(u',u'')=\cP(u',\langle u''l_{1},l_{1}\rangle,...,
\langle u''l_{m},l_{m}\rangle),
$$
where $\langle\cdot,\cdot\rangle$ is the scalar product in $\bR^{d}$.

The rest of the article is organized as follows.
In Section \ref{section 12.13.1} we show that
one may safely impose an additional assumption while
proving Theorem
\ref{theorem 9.23.1}. In Section \ref{section 12.13.2}
Theorem
\ref{theorem 9.23.1} is deduced from Theorem 
\ref{theorem 10.5.1} in which even more additional
assumptions are made. Then in Section
\ref{section 10.18.1}  the function $H$ is rewritten
in terms of pure second-order derivatives
along certain directions.

In a quite long Section \ref{section 9.22.1}
we consider finite-difference approximations
for equations   with ``constant''
coefficients and prove   interior estimates
for the   second-order differences
of solutions. In Section \ref{section 10.18.2}
we use the results of the previous section
in order to prove an analog of Theorem \ref{theorem 9.23.1}
for $H$, that include only pure second-order derivatives.
Here the reader will see the main underlying idea
of the paper, which roughly speaking
  is that on the set, say $\Gamma$,
 where the second-order
derivatives of $v$ are large we have $P[v]=K$
and   estimates similar to the ones
   from Section \ref{section 9.22.1}
show that the second order derivative on $\Gamma$
are controlled by their values on the boundary of $\Gamma$,
where they are under control by the definition of $\Gamma$.
Of course, the implementation of this idea requires
first proving that there are sufficiently regular
solutions of \eqref{9.23.2}. Since we do not know how to do 
that, we apply the above idea at the level of
finite-differences.

In the final short Section \ref{section 12.13.4}
we prove Theorem  \ref{theorem 10.5.1}.

In the proofs of various results in this article we use
the symbol $N$ sometimes with indices to denote constants
which may change from one occurrence to another and
we do not always specify on which data these  constants
depend. In these cases the reader should remember
that, if in the statement of a result there are constants
called $N$ which claimed to depend only on certain
parameters, then in the proof of the result
the constants $N$ also depend only on the same
parameters unless specifically stated otherwise.

\mysection{Reducing Theorem 
\protect\ref{theorem 9.23.1} to a particular case
where $-H_{u'_{0}}\geq \delta$}
                                   \label{section 12.13.1}

Suppose that Theorem \ref{theorem 9.23.1}
is true under the additional assumption that
\begin{equation}
                                                   \label{3.5.2}
-H_{u'_{0}}\geq \delta
\end{equation}
 at all points of
differentiability of $H(u,x)$ with respect to $u$.
 Then we are going to prove it
in the original form. Take
an $H$ satisfying only Assumption \ref{assumption 9.23.1},
take $n>0$,
and  consider the mapping $T_{n}:w\to v$ defined for any
$w\in C(\bar{\Omega})$ and mapping it into
a unique solution of
\begin{equation}
                                          \label{10.4.1}
\max(H[v]- v+n\chi(w/n),P[v]-K)=0
\end{equation}
in $\Omega$  (a.e.) with  boundary condition $v=g$,
where
$$
\chi(t)=(-1)\vee t\wedge 1.
$$
By assumption $v$ is well defined and $v=T_{n}w\in
C^{0,1}(\bar{\Omega})\cap C^{1,1}_{loc}(\Omega)$ and
$$
|v|,|D_{i}v|,\rho|D_{ij} v |\leq N(\bar{H}+n+K
+\|g\|_{C^{1,1}(\Omega)})\quad\text{in}
\quad \Omega \quad (a.e.), 
$$
$$
\|v\|_{W^{2}_{p}(\Omega)}\leq
N_{p}(\bar{H}+n+K+\|g\|_{W^{2}_{p}(\Omega)}).
$$
It follows that, for each $n$, $T_{n}$ maps $C(\bar{\Omega})$
into its compact subset.

\begin{lemma}
                                        \label{lemma 10.4.1}
For each $n$, the mapping $T_{n}$ is continuous
in $C(\bar{\Omega})$.
\end{lemma}

Proof. Let $w,w_{m}\in C(\bar{\Omega})$, $m=1,2,...$,
and assume that $\|w-w_{m}\|_{0,\Omega}\to0$ as $m\to\infty$,
where $\|\cdot\|_{0,\Omega}$ is the sup norm in
$C(\bar{\Omega})$. In light of uniqueness of solutions
of \eqref{10.4.1} with   boundary condition $v=g$,
to prove the lemma, it suffices to show that, at least
along a subsequence, $\|v-v_{m}\|_{0,\Omega}\to0$,
where $v=T_{n}w$, $v_{m}=T_{n}w_{m}$. Since $T_{n}
C(\bar{\Omega})$ is a compact set, there is a subsequence
and a $v \in C(\bar{\Omega})$
such that $\|v-v_{m}\|_{0,\Omega}\to0$ and $v=g$
on $\partial\Omega$. Without losing
generality we may assume that the above convergence holds
along the original sequence.
Now we need only  show  that $v=T_{n}w$.

Observe that for $m\geq r$ we have
$$
\max(H[v_{m}]- v_{m}+n\sup_{k\geq r}
\chi(w_{k}/n),P[v_{m}]-K)
\geq0 
$$
in $\Omega$ (a.e.).
Since the norms $\|v_{m}\|_{W^{2}_{d}(\Omega)}$
are bounded,
by Theorems 3.5.9 and 3.5.15 of \cite{Kr85}, whose
conditions are easily checked on the basis of
Remark \ref{remark 12.21.1},  we have  (a.e.)
$$ 
\max(H[v ]- v+n\sup_{k\geq r}\chi(w_{k}/n),P[v ]-K)
\geq0.
$$
By letting $r\to\infty$ we get  (a.e.)
$$ 
\max(H[v ]- v+n \chi(w /n),P[v ]-K)
\geq0.
$$
One obtains the opposite inequality  starting with
$$
\max(H[v_{m}]- v_{m}+n\inf_{k\geq r}\chi(w_{k}/n),P[v_{m}]-K)
\leq0.
$$
It follows that $v=T_{n}w$ indeed and the lemma is proved.

Now by Tikhonov's theorem we conclude that, for each $n$,
there exists $v^{n}\in C(\bar{\Omega})$ such that $v^{n}=
T_{n}v^{n}$.
By assumption $v^{n} \in
C^{0,1}(\bar{\Omega})\cap C^{1,1}_{loc}(\Omega)$ and
$$
|D_{i}v^{n}|,\rho|D_{ij} v^{n} |\leq N(\bar{H}+
\|v^{n}\|_{0,\Omega}
+K+\|g\|_{C^{1,1}(\Omega)})\quad\text{in}
\quad \Omega \quad (a.e.), 
$$
 \begin{equation}
                                           \label{10.4.3}
\|v^{n}\|_{W^{2}_{p}(\Omega)}\leq
N_{p}(\bar{H}+\|v^{n}\|_{0,\Omega} 
+K+\|g\|_{W^{2}_{p}(\Omega)}).
\end{equation}

\begin{lemma}
                                        \label{lemma 10.4.2}
There is a constant depending only on the diameter
of $\Omega$ and $\delta$ such that
$$
\|v^{n}\|_{0,\Omega}\leq N(\bar{H}+K+\|g\|_{C(\Omega)}).
$$
\end{lemma}

Proof. Introduce
$$
H^{n}_{K}(u,x)=\max(H(u,x)-u'_{0}+n\chi(u'_{0}/n),
P(u)-K)
$$
and observe that 
  $ H^{n}_{Ku'_{0}}\leq 0$ and
by  Hadamard's formula
$$
H^{n}_{K} (u' ,u'',x)-H^{n}_{K} (0,x)
= u''_{ij} \int_{0}^{1}H^{n} _{Ku''_{ij}}
(tu' ,tu'',x)\,dt
$$
$$
+\sum_{i\geq1}
 u' _{i} \int_{0}^{1}H^{n} _{Ku' _{i}}
(tu' ,tu'',x)\,dt+
 u' _{0} \int_{0}^{1}H^{n} _{Ku' _{0}}
(tu' ,tu'',x)\,dt.
$$
Then we see that, for each $n$,
 there exist $\cS_{\delta}$-valued
function $a $ and real-valued
 functions $b_{1},...,b_{d}$, $c$, and $f$
satisfying $|b_{i}|\leq \delta^{-1}$, 
$ c\geq0$, $|f|\leq \bar{H}+K$
such that in $\Omega$ (a.e.)
$$
a_{ij}D_{ij}v^{n}+b_{i}D_{i}v^{n}-cv^{n}=f.
$$
Now our result follows by the Alexandrov maximum
principle (see, for instance, Section 3.3 of \cite{Kr85}). 
The lemma is proved.

 Due to this lemma one can drop $\|v^{n}\|_{0,\Omega}$
in the right-hand sides of estimates \eqref{10.4.3}.
After that it only remains to observe that
for $n\geq\|v^{n}\|_{0,\Omega}$, the function
$v_{n}$ satisfies \eqref{9.23.2} since $\chi(v_{n}/n)=v_{n}/n$
and  Theorem 
\ref{theorem 9.23.1} holds in its original form.

Hence, in the rest of the article we suppose that
\eqref{3.5.2} holds at all points
of differentiability of $H$ with respect to $u$.

\mysection{Further reductions
of Theorem 
\protect\ref{theorem 9.23.1}}

                                   \label{section 12.13.2}

\noindent{\bf 1}. First, we show that we may additionally assume 
that for any
$x,y\in\bR^{d}$ and $u=(u',u'')$
\begin{equation}
                                                  \label{9.31.2}
|H(u,x)-H(u,y)|\leq N|x-y|(1+|u|),
\end{equation}
where $N$ is independent of $x,y,u$. 

Indeed, if Theorem \ref{theorem 9.23.1} is true
in this particular case, take a nonnegative
$\zeta\in C^{\infty}_{0}(\bR^{d})$, which integrates
to one,  set $\zeta^{n}(x)=n^{d}\zeta(nx)$, and introduce
$
H^{n}(u,x)
$ as the convolution of $H(u,x)$ and $\zeta^{n}$ performed
with respect to $x$. 
Observe that   $H^{n}$  satisfies
\eqref{3.5.2} and Assumption \ref{assumption 9.23.1}
with the same constant $\delta $, whereas
$$
|H^{n}(u,x)-H^{n}(u,y)|\leq  n |x-y|\sup_{z}|H(u,z)|
\sup| D\zeta|
$$
and \eqref{9.31.2} is satisfied since $ |H(u,z)|\leq
|H(0,z)|+N(d)\delta^{-1}|u|$. Then 
assuming that the assertions of Theorem \ref{theorem 9.23.1}
are true under our additional assumption,
we conclude that there exist solutions 
$v^{n} \in C^{0,1}(\bar{\Omega})\cap C^{1,1}_{loc}(\Omega)$
of
\begin{equation}
                                       \label{3.29.2}
\max(H^{n}[v^{n}],P[v^{n}]-K)=0
\end{equation}
in $\Omega$ (a.e.) with  boundary condition $v^{n}=g$,
for which estimates \eqref{1.13.1} and \eqref{1.13.2} hold with
$v^{n}$ in place of $v$ with the constants $N$ and $N_{p}$
from Theorem \ref{theorem 9.23.1} and with
$$
\overline{H^{n}}=\sup_{x\in\bR^{d}}|H^{n}(0,x)|
\quad\quad(\leq\bar{H})
$$
in place of $\bar{H}$. In particular,
\begin{equation}
                                              \label{10.3.4}
\check{H}^{n}_{K}[v^{m}]\geq0
\end{equation}
in $\Omega$ (a.e.) for all $m\geq n$, where
$$
\check{H}^{n}_{K}(u,x):=\sup_{k\geq n}
\max(H^{k}(u,x),P(u)-K).
$$

 Furthermore, being uniformly bounded and uniformly continuous,
the sequence $\{v^{n}\}$ has a subsequence uniformly converging
to a function $v$, for which \eqref{1.13.1} and \eqref{1.13.2}, 
of course, hold 
and $v  \in C^{0,1}(\bar{\Omega})\cap C^{1,1}_{loc}(\Omega)$.
In light of \eqref{10.3.4} and the fact that
 the norms $\|v^{n}\|_{W^{2}_{p}(\Omega)}$ are bounded,
by Theorems 3.5.9 and 3.5.15 of \cite{Kr85}
(the applicability of which is shown by an argument similar
to the one in Remark \ref{remark 2.28.1})
 we have 
\begin{equation}
                                              \label{10.3.2}
\check{H}^{n}_{K}[v ]\geq0 
\end{equation}
in $\Omega$ (a.e.).

  Then we
  notice that by the Lebesgue differentiation theorem
for any $u$
\begin{equation}
                                              \label{10.3.1}
\lim_{n\to\infty}\check{H}^{n}_{K}(u,x)=\max(H(u,x),P(u)-K)
\end{equation}
for almost all $x$. Since $\check{H}^{n}_{K}(u,x)$ are Lipschitz
continuous in $u$ with a constant independent of $x$ and $n$,
there exists a subset of $\Omega$ of full measure such that
\eqref{10.3.1} holds on this subset for all $u$.

We conclude that in $\Omega$ (a.e.)
\begin{equation}
                                              \label{10.3.3}
\max(H [v],P[v]-K)\geq0.
\end{equation}

The opposite inequality is obtained by considering
$$
\inf_{k\geq n}
\max(H^{k}(u,x),P(u)-K).
$$

\noindent{\bf  2}. Next, 
 we show that one may assume that $H$ is 
boundedly inhomogeneous with respect to $u$.
Introduce
$$
P_{0}(u)=\max_{a\in\cS_{\delta/2}}
\max_{\substack{|b_{i}|\leq2\delta^{-1}\\
i=1,...,d}}\max_{c\in[\delta/2,2
\delta^{-1}]}(a_{ij}u''_{ij}+b_{i}u'_{i}-cu'_{0}),
$$
where  the summations   are performed
 before the maximum is taken.
It is easy to see that $P_{0}[u]$ is a kind of Pucci's operator:
$$
P_{0}(u)=-(\delta/2)\sum_{k=1}^{d}\lambda_{k}^{-}(u'')
+2\delta^{-1}\sum_{k=1}^{d}\lambda_{k}^{+}(u'')
$$
$$
+2\delta^{-1}\sum_{k=1}^{d}|u'_{k}|-(\delta/2)(u'_{0})^{+}
+2\delta^{-1}(u'_{0})^{-},
$$
where $\lambda_{1}(u''),...,\lambda_{d}(u'')$ are 
the eigenvalues
of $u''$ and $a^{\pm}=(1/2)(|a|\pm a)$.

Recall that the function $P$ is introduced in the end of Section
\ref{section 2.5.1} and
observe that  
$$
 P(u)  =\max_{\substack{\hat{\delta}/2\leq a_{k}\leq
2\hat{\delta}^{-1} \\k=1,...,m} }
\max_{\substack{ |b_{i}|\leq
2\hat{\delta}^{-1}\\i=1,...,d } }
\max_{\hat{\delta}/ 2\leq c\leq2
\hat{\delta}^{-1}}\big[\sum_{i,j=1}^{d}
\sum_{k=1}^{m} a_{k}l_{ki}l_{kj}u''_{ij}
 +\sum_{i=1}^{d}b_{i}u'_{i} -cu'_{0}\big].
$$
Moreover,
 owing to property (ii) in the end of Section
\ref{section 2.5.1}, the collection of matrices
$$
\sum_{k=1}^{m}a_{k}l_{k}l_{k}^{*}
$$
such that $\hat{\delta}\leq a_{k}\leq
\hat{\delta}^{-1},k=1,...,m$, covers   $\cS_{\delta/4}$.
By combining this with  the fact that $\hat{\delta}
\leq\delta/2 $ (actually,
$\hat{\delta}
\leq\delta/4$, which will be used much later) we see that
$$
P(u)\geq -(\delta/4)\sum_{k=1}^{d}\lambda_{k}^{-}(u'')
+4\delta^{-1}\sum_{k=1}^{d}\lambda_{k}^{+}(u'')
$$
$$
+4\delta^{-1}\sum_{k=1}^{d}|u'_{k}|-(\delta/4)(u'_{0})^{+}
+4\delta^{-1}(u'_{0})^{-},
$$
\begin{equation}
                                          \label{3.6.1}
\geq P_{0}(u)+(\delta/4)\sum_{k=1}^{d}|\lambda_{k} (u'')|
+(\delta/4)\sum_{k=0}^{d}|u'_{k}|.
\end{equation}
   
In particular,
$P_{0}\leq P$ and
therefore,
$$
\max(H,P-K)=\max(H_{K},P-K),
$$
where $H_{K}=\max(H,P_{0}-K)$. It is easy to see that
the function $H_{K}$ satisfies Assumption
 \ref{assumption 9.23.1} and  \eqref{3.5.2}
 with $\delta/2$  in place of $\delta$.
It also satisfies \eqref{9.31.2} with the same constant $N$.

Furthermore, we have the following.   
\begin{lemma}
                                          \label{lemma 9.29.2}
There is a constant $\kappa>0$ depending only on $\delta$
and $d$ such that  for all $x\in\Omega$ and $u=(u',u'')$
\begin{equation}
                                            \label{9.29.2}
H \leq  P_{0}-\kappa
 \big(\sum_{i,j}|u''_{ij}| +\sum_{i}|u'_{i}|\big)+
 H(0,x),  
\end{equation}
 \begin{equation}
                                            \label{3.6.3}
H_{K} \leq  P -\kappa
 \big(\sum_{i,j}|u''_{ij}| +\sum_{i}|u'_{i}|\big)+
\ H^{+}(0,x).
\end{equation}
Furthermore,
$$
 H(u,x) \leq N\big(\sum_{i,j}|u''_{ij}| +\sum_{i}|u'_{i}|\big)+
 H(0,x),
$$
$$
|H(u,x)|\leq N\big(\sum_{i,j}|u''_{ij}| +\sum_{i}|u'_{i}|\big)+
|H(0,x)|,
$$
where the constant $N$ depends only on  $\delta$.
\end{lemma}

Proof.  Observe that if a number $p\in( a,b)$, $a<b$, and 
$y\in\bR$, then
$$
yp\leq y^{+}b-y^{-}a.
$$
Then from  Hadamard's formula
$$
H (u' ,u'',x)-H (0,0,x)
=u''_{ij}\int_{0}^{1}H _{u''_{ij}}
(tu' ,tu'',x)\,dt
$$
$$
+\sum_{i\geq1}
u' _{i}\int_{0}^{1}H _{u' _{i}}
(tu' ,tu'',x)\,dt+
u' _{0}\int_{0}^{1}H _{u' _{0}}
(tu' ,tu'',x)\,dt
$$
we obtain
$$
H (u' ,u'',x)-H (0,0,x)\leq
\delta^{-1}\sum_{k}\lambda^{+}_{k}(u'')-
\delta \sum_{k}\lambda^{-}_{k}(u'') 
$$
$$
+\delta^{-1}\sum_{i\geq1}|u' _{i}|
-\delta (u'_{0})^{+}+\delta^{-1}
(u'_{0})^{-}=P _{0}(u' ,u'')
$$
$$
-\delta^{-1}\sum_{k}\lambda^{+}_{k}(u'')
-(\delta /2)\sum_{k}\lambda^{-}_{k}(u'') 
-\delta^{-1}\sum_{i\geq1}|u'_{k}|-
\delta^{-1}
(u'_{0})^{-}-(\delta/2) (u'_{0})^{+}
$$
and \eqref{9.29.2} follows since
$$
\big[\sum_{k}(\lambda^{+}_{k}(u'') +\lambda^{-}_{k}(u'') )
\big]^{2}=
\big(\sum_{k}|\lambda_{k}(u'')| \big)^{2}
$$
$$
\geq\sum_{k}|\lambda_{k}(u'')|^{2}=
\sum_{i,j}|u''_{ij}|^{2}\geq d^{-2}
\big(\sum_{i,j}|u''_{ij}|\big)^{2}.
$$

Estimate \eqref{3.6.3} follows from \eqref{9.29.2} and
\eqref{3.6.1}. Finally,
the second assertion of the lemma
follows directly from the above Hadamard's formula.
 The lemma is proved.

In addition, $H_{K}$ is boundedly inhomogeneous 
with respect to $u$ in the sense
that at all points of differentiability of 
$H_{K}(u,x)$ with respect to $u$
\begin{equation}
                                              \label{9.30.01}
|H_{K}(u,x)-H_{Ku''_{ij}}(u,x)u''_{ij}
-H_{Ku'_{r}}(u,x)u'_{r}|\leq N(|H_{K}(0,x)| +K),
\end{equation}
where $N$ depends only on $\delta$ and $d$.

Indeed, if
\begin{equation}
                                            \label{9.30.2}
\kappa
 \big(\sum_{i,j}|u''_{ij}| +\sum_{i}|u'_{i}|\big) 
 \geq  H^{+}(0,x)+K,
\end{equation}
then by Lemma \ref{lemma 9.29.2} 
$$
H (u,x)\leq P_{0}(u)-\kappa
 \big(\sum_{i,j}|u''_{ij}| +\sum_{i}|u'_{i}|\big)+
 H^{+}(0,x)\leq P_{0}(u)-K,
$$
so that
$H_{K}(u,x)=P_{0}(u)-K$  
and the
left-hand side of \eqref{9.30.01} is just $K$
owing to  the  fact that $P_{0}$ is positive homogeneous
of degree one. On the other hand, if the opposite inequality
holds in \eqref{9.30.2}, then  again in light of Lemma
\ref{lemma 9.29.2}  
the
left-hand side of \eqref{9.30.01} is dominated by
$$
N\big(\sum_{i,j}|u''_{ij}| +\sum_{i}|u'_{i}|\big)+
|H_{K}(0,x)|\leq N(|H_{K}(0,x)|+ H^{+}(0,x)+K),
$$
where 
$$
H(0,x)\leq\max(H(0,x),-K)=H_{K}(0,x),\quad
H^{+}(0,x)\leq |H_{K}(0,x)|.
$$

Furthermore, as we have noticed above $H_{K}$ satisfies
Assumption \ref{assumption 9.23.1}  
and \eqref{3.5.2} (with $\delta/2$ in place of $\delta$)
and as is easy to see $|H_{K}[0]|\leq |H[0]|+K$,
which shows that in the rest of the article we may (and will) assume
that not only Assumption \ref{assumption 9.23.1} and
 \eqref{3.5.2} are satisfied with $\delta/2$ in place of $\delta$
and \eqref{9.31.2} holds with a constant $N$,
but also at all points of differentiability of $H$
with respect to $u$
\begin{equation}
                                              \label{3.5.1}
|H (u,x)-H_{ u''_{ij}}(u,x)u''_{ij}
-H_{ u'_{r}}(u,x)u'_{r}|\leq N_{0} ,
\end{equation}
where $N_{0}$ is a constant and
 \begin{equation}
                                            \label{3.6.4}
H  \leq  P -\kappa
 \big(\sum_{i,j}|u''_{ij}| +\sum_{i}|u'_{i}|\big)+
|H (0,\cdot)|,
\end{equation}
where $\kappa$ is the constant from Lemma \ref{lemma 9.29.2}.
  By the way we keep track
of the value of $\delta$ in Assumption \ref{assumption 9.23.1}  
and \eqref{3.5.2} because $P(u)$ is already fixed and defined
by $d$ and $\delta$.

\noindent{\bf 3}. Finally, we show that we may assume that, 
\begin{equation}
                                               \label{10.5.2}
H(u,x)={\rm tr}\,u''-u_{0}'
\,\,\text{for all $u$ if $x$ is in a neighborhood of}\quad
\partial \Omega,
\end{equation}
that is,
for an  
$\varepsilon>0$, we have $H(u,x)={\rm tr}\,
u''-u'_{0}$ if $\rho(x)\leq\varepsilon$. Indeed,
take a continuous function $\zeta(t)$, $t\geq0$ such that
$\zeta(t)=1$ for $t\in[0,1]$, $\zeta(t)=0$
for $t\geq 2$, and $0\leq\zeta\leq1$. Introduce
$$
H^{1/\varepsilon}(u,x)=(1-\zeta(\rho(x)/\varepsilon))
H (u,x)+\zeta(\rho(x)/\varepsilon)({\rm tr}\,
u''-u'_{0}).
$$

Notice that  $H^{1/\varepsilon}$ satisfies Assumption
\ref{assumption 9.23.1} and
 \eqref{3.5.2}  with $\delta/2$ in place of $\delta$,
satisfies \eqref{9.31.2} with a constant $N$ depending on
$\varepsilon$ but independent of $x,y,u$,
and satisfies \eqref{3.5.1} with the same constant $N_{0}$.
As long as \eqref{3.6.4} is concerned,
 observe that by Lemma \ref{lemma 9.29.2}
applied to $H={\rm tr}\,
u''-u'_{0}$ and by the inequality $P_{0}\leq P$ we have
$$
{\rm tr}\,
u''-u'_{0}\leq P -\kappa
 \big(\sum_{i,j}|u''_{ij}| +\sum_{i}|u'_{i}|\big).
$$
Then owing to \eqref{3.6.4}
$$
H^{1/\varepsilon}\leq P -\kappa
 \big(\sum_{i,j}|u''_{ij}| +\sum_{i}|u'_{i}|\big)+
|(1-\zeta(\rho /\varepsilon))H(0,\cdot)|
$$
$$
=P -\kappa
 \big(\sum_{i,j}|u''_{ij}| +\sum_{i}|u'_{i}|\big)+
|H^{1/\varepsilon}(0,\cdot)|.
$$

Therefore, if the assertions of Theorem \ref{theorem 9.23.1} 
hold  under the above additional assumptions,  then we have
a sequence of functions $v^{n} \in C^{0,1}(\bar{\Omega})
\cap C^{1,1}_{loc}(\Omega)$ satisfying \eqref{3.29.2}
(with new $H^{n}=H^{1/\varepsilon}$ for $\varepsilon=1/n$).

After  that by repeating literally the argument  in 
no.~{\bf 1} we come to \eqref{10.3.2} and 
since, obviously, $H^{1/\varepsilon}(u,x)\to H(u,x)$
as $\varepsilon\downarrow0$ for any $x\in\Omega$, 
we conclude that equation \eqref{10.3.3} holds (a.e.)
 and we finish the argument
as in no.~{\bf 1}. 
 
As a result of the above arguments we see that to 
prove Theorem \ref{theorem 9.23.1} it suffices to prove
the following.
\begin{theorem}
                                   \label{theorem 10.5.1}
Suppose that Assumption \ref{assumption 9.23.1}
is satisfied with $\delta/2$ in place of
$\delta$. Also assume that  \eqref{3.6.4} and \eqref{10.5.2} 
 hold. Finally, assume that
estimate \eqref{9.31.2} holds
for any
$x,y\in\bR^{d}$ and $u=(u',u'')$ with 
a constant $N$
 and \eqref{3.5.2} and \eqref{3.5.1}
hold  at all points of differentiability of   
$H (u,x)$ with respect to $u$.

Then the assertions of Theorem \ref{theorem 9.23.1}  
hold true with $P$ introduced in the end of
Section \ref{section 2.5.1}.
\end{theorem}

\mysection{Writing $H$ in Theorem \protect\ref{theorem 10.5.1}
in a special form}

                                        \label{section 10.18.1} 
Here we suppose that the assumptions
of Theorem \ref{theorem 10.5.1} are satisfied
and take the objects introduced in the end
of Section \ref{section 2.5.1}. Owing to the
the assumptions
of Theorem \ref{theorem 10.5.1} by Theorem   7.1 of \cite{Kr11}
(see the beginning of its proof in \cite{Kr11})
 there exists a
function
$\cH(z,x)$ defined for  
$$
z  = ( z ',z ''),\quad z '=
(z '_{0},...,z '_{d})
\in\bR^{d+1},\quad z ''\in\bR^{m} ,\quad x\in\bR^{d}
$$
such that: 

(i) The function $\cH $ is Lipschitz continuous in $z $
with Lipschitz constant $\hat{\delta}^{-1}$ and there exists
a constant $N'$ such that
$$
|\cH(z,x)-\cH(z,y)|\leq N'|x-y|(1+|z|) 
$$
for all $x,y\in\bR^{d}$ and $z$.

(ii) We have $\cH(z,x )=H(u,x)$ if 
$z '=u'$ and for all $j=1,...,m$
$$
z ''_{j}=\langle u''l_{j},l_{j}\rangle.
$$
In particular, $\cH(0,x)=H(0,x)$ and
if $v(x)$ is a real-valued function
which is twice differentiable
at a point $x\in\bR^{d}$, at this point we have
$$
H[v](x)=\cH[v](x)
$$
where
$$
\cH[v](x) =\cH(v,Dv,D^{2}_{l_{1}}v,...,D^{2}_{l_{m}}v,x),
\quad D^{2}_{l}v=v_{x_{i}x_{j}}l_{i}l_{j}.
$$

(iii) At all points $(z,x) $ at which $\cH (z,x)$ is differentiable
with respect to $z$
we have 
\begin{equation}
                                                  \label{5.10.1}
|\cH _{z ' _{i}}(z,x )|\leq 4\delta ^{-1},\quad i=1,...,d,
\end{equation}
\begin{equation}
                                                  \label{5.10.2} 
 \delta/4\leq-\cH_{z ' _{0}}(z,x )\leq 4\delta ^{-1},
\quad
\hat{\delta}^{-1}\geq \cH_{z ''_{j}}(z,x )\geq\hat{\delta},\quad
j=1,...,m .
\end{equation}

The proofs in \cite{Kr11}
use the fact that \eqref{3.5.1} holds and yield
the function $\cH$ such that, in addition, at all points 
$(z,x) $ at which $\cH (z,x)$ is differentiable
with respect to $z$
we also have 
$$
 |\cH(z,x )-\langle z , D_{z} \cH
 (z,x )\rangle |\leq 2N_{0}.
$$
However, the latter property of $\cH$ will not be used 
in the future, so that we only used assumption
\eqref{3.5.1} to be sure that $\cH$ with the
properties (i)-(iii) exists.

\mysection{An auxiliary equation}

                                       \label{section 9.22.1}
Some notation in this section
are different from the previous ones.
Fix an $h\in(0,1]$ and for $\xi\in\bR^{d}$
and any function $\phi$ on $\bR^{d}$ introduce
$$
T_{\xi}\phi(x)=\phi(x+h\xi),\quad
\delta_{\xi}=h^{-1}(T_{\xi}-1),\quad\Delta_{\xi}
=h^{-2}(T_{\xi}-2+T_{-\xi}).
$$
Notice that $h$ enters the definition of $T_{\xi}$
and $\delta_{\xi}$ and $\Delta_{\xi}$ are usual
approximations for the first and second-order derivative
along $\xi$.

Let $m\geq1$ be an integer
and let $\ell_{-m},...,\ell_{-1},\ell_{1},...,\ell_{m}$
be some fixed vectors in $\bR^{d}$
such that 
$$
\ell_{-k}=-\ell_{k}.
$$ 
Next denote $\Lambda=\{\ell_{k}:k=\pm1,...,\pm m\}$,
$$
\Lambda_{1}= \Lambda,\quad \Lambda_{n+1}= 
 \Lambda_{n}+ \Lambda ,\quad n\geq1,
\quad \Lambda_{\infty}=\bigcup_{n}\Lambda_{n}\,,
$$

Let $m'\geq0$ be an integer $\leq m$ and let $A=\{\alpha=(a,b,c)\}$
be a closed bounded set in $\bR^{2m }\times\bR^{m'}\times
\bR$, so that 
$$
a=(a_{-m},a_{-m+1},...,a_{-1},a_{1},...,a_{m})\in\bR^{2m},
$$
$$
b=
( b_{1},...,b_{m'})\in\bR^{ m'},
$$ 
and $c\in\bR$.
Also let
$f(\alpha,x)$ be
a real-valued function defined for $\alpha\in A$, $x\in\bR^{d}$.

Fix an $r\in\{1,...,m\}$ and for $k=\pm 1,...,\pm m$ set
$$
\delta_{h,k}=\delta_{k}=\delta_{\ell_{k}},\quad 
\Delta_{h,k}=\Delta_{k}=\Delta_{\ell_{k}}.
$$ 

\begin{assumption}
                                   \label{assumption 9.18.1}
There are constants $\delta>0$ and $K_{1},K_{2}
\in[0,\infty)$ 
  such that

(i) For any $(a,b,c)\in A$ and all $k$ we have
$$
a_{k}=a_{-k},\quad
\delta \leq a_{k}\leq\delta^{-1},\quad|b_{k}|\leq\delta^{-1},
\quad b_{k}^{-}\leq ha_{k},
\quad c\geq0;
$$

(ii) The function $f$ is continuous in 
$\alpha$ for any $x$ and $|\delta_{r}f|\leq K_{1}$,
$\Delta_{r}f\geq-K_{2}$ on $\bR^{d}$.

\end{assumption}

For $u=(u',u'')$ with
$$
u'=(u'_{0}, u'_{ 1},...,u'_{m'}),\quad
 u'' =(u''_{-m},...,u''_{-1},u''_{1},...,u''_{m}),
$$
introduce
$$
\cP(u,x)=\max_{\alpha=(a,b,c)\in A}\big(\sum_{|k|=1}^{m} 
a_{k}u''_{k}+\sum_{ k =1}^{m'}
b_{k}u'_{k} -cu'_{0}+f(\alpha,x)\big).
$$

For any function $u$ on $\bR^{d}$ define
$$
\cP[u](x)=\cP(u(x),\delta u(x),\delta^{2}u(x) ,x),
$$
where
$$
\delta u=( \delta_{1}u ,...,
\delta_{m'}u ),
$$
$$
\delta^{2} u=(\Delta_{-m}u ,...,
\Delta_{-1}u ,\Delta_{1}u ,...,
\Delta_{m}u).
$$
In connection with this notation a natural question arises as
to why use $\ell_{k}$ along with $\ell_{-k}=-\ell_{k}$ since
$\Delta_{k}=\Delta_{-k}$ and 
$$
a_{k}\Delta_{k}=2\sum_{k\geq1}a_{ k}\Delta_{k}
$$
owing to the assumption that $a_{k}=a_{-k}$. This is done
for the sake of convenience of computations. For instance,
$$
\Delta_{k}(uv)=u\Delta_{k}v+v\Delta_{k}+(\delta_{k}u)(\delta_{k}v)
+(\delta_{-k}u)(\delta_{-k}v)
$$
(no summation in $k$). At the same time
$$
a_{k}\Delta_{k}(uv)=ua_{k}\Delta_{k}v+va_{k}\Delta_{k}
+2a_{k}(\delta_{k}u)(\delta_{k}v)
$$
as if we were dealing with usual partial derivatives.

Fix a nonempty finite set $Q^{o}\in h\Lambda_{\infty}$ and 
let
$$
Q := Q^{o}\cup\{x +h\Lambda   :x\in Q^{o}\}.
$$
Next take a function $\eta\in C^{\infty}(\bR^{d})$
with bounded derivatives,
such that $|\eta|\leq1$ and set $\zeta=\eta^{2}$,
$$
|\eta'(x)| =|\eta'(x)|_{h} =\sup_{  k}|\delta_{ k}\eta  (x)|,\quad
|\eta''(x)| =|\eta''(x)|_{h} =\sup_{ k}| 
\Delta_{ k}\eta  (x)|,
$$
$$
\|\eta'\|=\|\eta'\|_{h}=\sup_{h\Lambda_{\infty}}|\eta' |_{h},\quad
\|\eta''\|=\|\eta''\|_{h}=\sup_{h\Lambda_{\infty}}|\eta'' |_{h},
$$

Finally, let $u$ be a function on $\bR^{d}$
such that in $Q^{o}$
\begin{equation}
                                                \label{9.19.9}
\cP[u]=0
\end{equation}
 and $\cP[u]\leq0$ on $Q\setminus Q^{o}$.

\begin{theorem}
                                      \label{theorem 9.18.1}
There  exist constants $N=N(m,\delta)\geq1$
 and $N^{*}=N^{*}(m,\delta)$
such that for   any constant $\nu$
satisfying
$$
\nu\geq N^{*} \|\eta'\|+N (\|\eta''\|+\|\eta'\|^{2}),
$$
we have in $Q$ that  (recall that  $a^{\pm}=(1/2)(|a|\pm a)$) 
\begin{equation}
                                               \label{9.18.3}
\zeta^{2}[ ( \Delta_{r}u)^{-}]^{2}\leq
\max_{Q\setminus Q^{o}}\zeta^{2}[ ( \Delta_{r}u)^{-}]^{2}
+(N\nu+N^{*})\bar{W}_{r}+N\nu^{-2}K_{2}^{2} 
+ \nu^{-1}K_{1}^{2},
\end{equation}
where
$$
\bar{W}_{r}=\max_{Q}(|\delta_{r}u|^{2}+
|\delta_{-r}u|^{2}).
$$
Furthermore, $N^{*}=0$ if $b\equiv0$.  

\end{theorem}

\begin{remark}
Theorem \ref{theorem 9.18.1} looks very much like
Theorem 1.1 of \cite{Kr2}. However, in the latter the boundary
of $Q^{o}$ is ``twice fatter'' and   
all mixed second-order differences are present under 
the maximum 
sign in the corresponding counterpart of
\eqref{9.18.3}. Our idea is to apply
 Theorem \ref{theorem 9.18.1} to regions where
at least one of pure second-order differences is
large. Then outside  the region all of them
will be under control. Yet this does not provide
any control of {\em mixed\/} differences on the boundary
of the region 
and makes it impossible to apply Theorem 1.1 of \cite{Kr2},
where the driving goal was to obtain estimates
for equations with variable coefficients
and estimating all mixed second-order finite differences
was necessary.

\end{remark}

In the following arguments no
 summation with respect to $r$ is done. The number $r$
 is 
fixed in the beginning of the section.
For simplicity of notation set
$$
u_{rr}=\Delta_{r}u,\quad
u_{r}=\delta_{r}u,\quad u_{kr}
=-\delta_{-k} \delta_{r}u.
$$
Notice that in the above line
the last notation when $k=r$ is consistent with the first
one.

In the following two lemmas the fact that $u$ is a solution
of \eqref{9.19.9} is not used and 
$$
u^{-}_{rr}=
(u_{rr})^{-}.
$$
\begin{lemma}
                                        \label{lemma 9.19.1}
There exists $N=N(m,\delta)$ and $N^{*}=N^{*}(m,\delta)$
such that, if $N^{*}h\leq 1$, on $Q^{o}$
for any $\alpha=(a,b,c)\in A$ we have
$$
- 2\zeta^{2}u_{rr}^{-}[a_{k}\Delta_{k}+b_{k}\delta_{k}]u_{rr}^{-}
\geq-[a_{k}\Delta_{k}+b_{k}\delta_{k}]
(\zeta^{2}(u_{rr}^{-})^{2}) 
$$

\begin{equation}
                                                \label{9.19.8}
-N^{*}|\eta'|\zeta a_{k}u_{kr}^{2}
-N(|\eta''|+ |\eta'|^{2})\zeta (u_{rr}^{-})^{2} 
-(N ^{*}|\eta'|^{2}+N|\eta'|^{4})\bar{W}_{r},
\end{equation}

$$
-2\zeta  u_{r}[a_{k}\Delta_{k}+b_{k}\delta_{k}]u_{r}
\geq -[a_{k}\Delta_{k}+b_{k}\delta_{k}](\zeta u_{r}^{2})
$$

\begin{equation}
                                                \label{9.19.7}
+ \zeta a_{k}u_{kr}^{2}- N(|\eta''|+|\eta'|^{2}
 )\bar{W}_{r}-N^{*}|\eta'|\bar{W}_{r}.
\end{equation}
Furthermore, $N^{*}=0$ if $b\equiv0$.

\end{lemma}

Proof.
 Here is the result  of simple computations,
which can be found, for instance, in the proof of Lemma
5.1 of \cite{Kr2}. For any $\alpha\in A$ we have
$$
a_{k}\Delta_{k}(\zeta^{2}( u_{rr }^{-})^{2})=
2\zeta^{2}u_{rr }^{-}a_{k}\Delta_{k}u_{rr }^{-}
+2a_{k}[\delta_{k}(\zeta u_{rr}^{-})]^{2}
$$

\begin{equation}
                                                 \label{9.18.5}
+4a_{k}u_{rr }^{-}(\delta_{k}\zeta)\delta_{k}(\zeta u_{rr }^{- })
+2a_{k}( u_{rr }^{-})^{2}[\zeta\Delta_{k}\zeta
-2(\delta_{k}\zeta)^{2}]
-4hu_{rr }^{-}a_{k}(\delta_{k}\zeta)^{2}
\delta_{k}u_{rr }^{-}.
\end{equation}
We also know from Lemma 4.3 of \cite{Kr2} that

$$
|\Delta_{k}\zeta|\leq2(|\eta''|+|\eta'|^{2}),
\quad (\delta_{k}\zeta)^{2}
\leq N|\eta'|^{2}\zeta+Nh^{2}|\eta'|^{4}\leq N |\eta' |^{2}.
$$

It follows that, for any $\varepsilon>0$,

$$
|4a_{k}u_{rr }^{-}(\delta_{k}\zeta)\delta_{k}(\zeta u_{rr }^{- })|
\leq \varepsilon a_{k} [\delta_{k}(\zeta u_{rr}^{-})]^{2}
+N\varepsilon^{-1}(u_{rr}^{-})^{2}( |\eta'|^{2}
\zeta+ h^{2}|\eta'|^{4}),
$$
where
 
\begin{equation}
                                                 \label{11.17.9}
(u_{rr}^{-})^{2}h^{2}\leq|h u_{rr}|^{2}=|(T_{\ell_{r}}-1)
u_{-r}|^{2},
\end{equation}
so that $|\eta'|^{4}(u_{rr}^{-})^{2}h^{2}\leq 4
|\eta'|^{4}\bar{W}_{r}$ in $Q^{o}$. Therefore,
in $Q^{o}$

 $$
|4a_{k}u_{rr }^{-}(\delta_{k}\zeta)\delta_{k}(\zeta u_{rr }^{- })|
\leq \varepsilon a_{k} [\delta_{k}(\zeta u_{rr}^{-})]^{2}
+N\varepsilon^{-1}|\eta'|^{2}\zeta (u_{rr}^{-})^{2} 
+N\varepsilon^{-1}\bar{W}_{r}|\eta'|^{4} .
$$

Similarly

$$
|2a_{k}( u_{rr }^{-})^{2}[\zeta\Delta_{k}\zeta
-2(\delta_{k}\zeta)^{2}]|\leq
N(|\eta''|+|\eta'|^{2})\zeta (u_{rr }^{-})^{2}
+N \bar{W}_{r}|\eta'|^{4} .
$$

By Lemma 4.3 of \cite{Kr2} for any $\varepsilon\in(0,1]$

$$
h(\delta_{k}\zeta)^{2}|u_{rr}^{-}\delta_{i}u_{rr}^{-}|\leq
\varepsilon|\delta_{i}(\zeta u_{rr}^{-})|^{2}+\varepsilon|\eta'|^{2}
\zeta(u_{rr}^{-})^{2}
$$

$$
+N\varepsilon^{-1}
|\eta'|^{4}[(hu_{rr}^{-})^{2}+(h^{2}\delta_{i}
u_{rr}^{-})^{2} ].
$$

Estimate \eqref{11.17.9} leads to

$$
h(\delta_{k}\zeta)^{2}|u_{rr}^{-}\delta_{i}u_{rr}^{-}|\leq
\varepsilon|\delta_{i}(\zeta u_{rr}^{-})|^{2}+ |\eta'|^{2}
\zeta(u_{rr}^{-})^{2}
$$

$$
+N\varepsilon^{-1}
\bar{W}_{r}|\eta'|^{4}+N\varepsilon^{-1}|\eta'|^{4}
 (h^{2}\delta_{i}
u_{rr}^{-})^{2}
$$
on $Q^{o}$, where the last term is estimated by using the
fact that $|\delta_{i}\psi^{-}|\leq|\delta_{i}\psi|$ 
for any function $\psi$ implying that

$$
|\eta'|^{4}(h^{2}\delta_{i}
u_{rr}^{-})^{2}) \leq|\eta'|^{4}|h\delta_{i}(u_{r}+u_{-r})|^{2}
=|\eta'|^{4}|(T_{\ell_{i}}-1)(u_{r}+u_{-r})|^{2}\leq N
|\eta'|^{4}\bar{W}_{r}.
$$
Hence on $Q^{o}$

$$
h(\delta_{k}\zeta)^{2}|u_{rr}^{-}\delta_{i}u_{rr}^{-}|\leq
\varepsilon|\delta_{i}(\zeta u_{rr}^{-})|^{2}+ |\eta'|^{2}
\zeta(u_{rr}^{-})^{2}
+N\varepsilon^{-1}
\bar{W}_{r}|\eta'|^{4}.
$$
Upon combining these estimates, choosing $\varepsilon$
appropriately, and coming back to
\eqref{9.18.5}, 
we find on $Q^{o}$ that

$$
 - 2\zeta^{2}u_{rr}^{-}a_{k}\Delta_{k}u_{rr}^{-}
\geq-a_{k}\Delta_{k}(\zeta^{2}(u_{rr}^{-})^{2})
+a_{k}[\delta_{k}(\zeta u_{rr}^{-})]^{2}
$$

\begin{equation}
                                                 \label{9.18.8}
-N(|\eta''|+ |\eta'|^{2})\zeta (u_{rr}^{-})^{2} 
-N \bar{W}_{r}|\eta'|^{4}.
\end{equation}
Next,

$$
b_{k}\delta_{k}(\zeta^{2}(u_{rr}^{-})^{2}))
=2\zeta u_{rr}^{-}b_{k}\delta_{k}(\zeta u_{rr}^{-})
+b_{k}h[\delta_{k}(\zeta u_{rr}^{-})]^{2}
$$

$$
=2\zeta^{2}u_{rr}^{-}b_{k}\delta_{k}u_{rr}^{-}
+2\zeta(u_{rr}^{-})^{2}b_{k}\delta_{k}\zeta
+2hb_{k}u_{rr}^{-}\zeta(\delta_{k}\zeta)\delta_{k}u_{rr}^{-}
+b_{k}h[\delta_{k}(\zeta u_{rr}^{-})]^{2}.
$$

Here $|\delta_{k}\zeta|\leq 2|\eta'|$, since $|\eta|\leq 1$.
Also  $a_{k}\geq\delta$, so that

$$
|2\zeta(u_{rr}^{-})^{2}b_{k}\delta_{k}\zeta|
\leq N^{*}|\eta'|\zeta a_{k}u_{kr}^{2}.
$$

Furthermore, $hT_{\ell_{k}}u^{-}_{rr}
=T_{\ell_{k}}(u_{r}+u_{-r})^{-}$, implying that on $Q^{o}$

$$
|2hb_{k}u_{rr}^{-}\zeta(\delta_{k}\zeta)\delta_{k}u_{rr}^{-}|
=|2hb_{k}u_{rr}^{-} (\delta_{k}\zeta)[\delta_{k}(\zeta u_{rr}^{-})
-(\delta_{k}\zeta)T_{\ell_{k}}u_{rr}^{-}]|
$$

$$
\leq N^{*}|\eta'|\bar{W}_{r}^{1/2}a_{k}^{1/2}
|\delta_{k}(\zeta u_{rr}^{-})|+N^{*}
u_{rr}^{-}\sum_{k}(\delta_{k}\zeta
)^{2}\bar{W}_{r}^{1/2},
$$
where, owing to the inequality $h|\eta'|\leq2$,
 the last term is dominated by

$$
N^{*} \bar{W}_{r}^{1/2}u_{rr}^{-}
(\zeta|\eta'|^{2}+h^{2}|\eta'|^{4})
\leq (N^{*} \zeta^{1/2}|\eta'|^{3/2}\bar{W}_{r}^{1/2})
(|\eta'|^{1/2}\zeta^{1/2}u_{rr}^{-})
$$

$$
+N^{*} |\eta'|^{3}\bar{W}_{r}
\leq N^{*}|\eta'|\zeta a_{k}u_{kr}^{2}+N^{*} |\eta'|^{3}\bar{W}_{r}.
$$

Hence,   in $Q^{o}$

$$
-2\zeta^{2}u_{rr}^{-}b_{k}\delta_{k}u_{rr}^{-}
\geq  -b_{k}\delta_{k}(\zeta^{2}(u_{rr}^{-})^{2}))
-N^{*}|\eta'|\zeta a_{k}u_{kr}^{2}
$$

\begin{equation}
                                                 \label{9.19.2}
-(N^{*}h+1/2)a_{k}[\delta_{k}(\zeta u^{-}_{rr})]^{2}-N^{*} 
(|\eta'|^{2}+|\eta'|^{3})\bar{W}_{r}.
\end{equation}
For $N^{*}h\leq1/2$ estimates \eqref{9.18.8} and \eqref{9.19.2}
and the fact that $|\eta'|^{3}\leq|\eta'|^{2}+|\eta'|^{4}$
lead to \eqref{9.19.8}.

To prove \eqref{9.19.7} observe that
(recall that $\delta_{k}u_{r}=-u_{-kr}$
and $a_{-k}=a_{k}$)

$$
a_{k}\Delta_{k}(\zeta u_{r}^{2})=a_{k}\zeta[2u_{r}\Delta_{k}u_{r}
+2u_{kr}^{2}]+a_{k}u_{r}^{2}\Delta_{k}\zeta
+2a_{k}(\delta_{k}\zeta)( hu_{-kr}^{2}-2u_{r}u_{-kr}),
$$
where

$$
u_{r}^{2}|\Delta_{k}\zeta|\leq N(|\eta''|+|\eta'|^{2})\bar{W}_{r},
$$
and owing to \eqref{11.17.9} 

$$
|2a_{k}(\delta_{k}\zeta)( hu_{-kr}^{2}-2u_{r}u_{-kr})|
\leq Na_{k}(|\eta'|\zeta^{1/2}+h|\eta'|^{2})
(|u_{r}u_{kr}|+h u_{kr}^{2}) 
$$

$$
\leq N a_{k}|u_{kr}|(|\eta'|\zeta^{1/2}+h|\eta'|^{2})
\bar{W}_{r}^{1/2}
\leq N a_{k}|u_{kr}| |\eta'|\zeta^{1/2}\bar{W}_{r}^{1/2}
+ N|\eta'|^{2}\bar{W}_{r}
$$

$$
\leq (1/2)\zeta a_{k}u_{kr}^{2}+ N|\eta'|^{2}\bar{W}_{r}.
$$

It follows that in $Q^{o}$

\begin{equation}
                                                 \label{9.19.1}
-2\zeta a_{k} u_{r}\Delta_{k}u_{r}
\geq -a_{k}\Delta_{k}(\zeta u_{r}^{2})
+(3/2)\zeta a_{k}u_{kr}^{2}- N(|\eta''|+|\eta'|^{2})\bar{W}_{r}.
\end{equation}
 
Also

$$
b_{k}\delta_{k}(\zeta u_{r}^{2})=
 2\zeta u_{r} b_{ k}\delta_{k} u_{r}+\zeta b_{-k}h u_{kr}^{2}
+u_{r}^{2}b_{k}\delta_{k}\zeta
+hb_{k}(\delta_{k}\zeta)[ hu_{-kr}^{2}-2u_{r}u_{-kr}].
$$
Here in $Q^{o}$

$$
|u_{r}^{2}b_{k}\delta_{k}\zeta|\leq N^{*}|\eta'|W_{r},
\quad
|hb_{k}(\delta_{k}\zeta)[ hu_{-kr}^{2}-2u_{r}u_{-kr}]
\leq N^{*}|\eta'|\bar{W}_{r},
$$
where the  last estimate follows from an equality similar to
\eqref{11.17.9}. Furthermore,

$$
|\zeta b_{-k}h u_{ kr}^{2}|\leq (1/2)\zeta a_{k}u_{kr}^{2}
$$
if $N^{*}h\leq 1$ and $N^{*}$ is chosen
appropriately. 

Upon combining this estimates with \eqref{9.19.1} we come to 
\eqref{9.19.7} and the lemma is proved.

For a constant $\nu\geq0$ introduce an operator
(recall that $r$ is fixed)

$$
L_{\nu}\phi=\zeta^{2}u_{rr}^{-}\Delta_{r}\phi
-\nu\zeta u_{r}\delta_{r}\phi.
$$
Observe that
\begin{equation}
                                              \label{9.21.1}
L_{\nu}u=-\zeta^{2}(u_{rr}^{-})^{2}-\nu\zeta u_{r}^{2}
=:-V_{\nu}.
\end{equation}

\begin{lemma}
                                        \label{lemma 9.20.1}
There exists $N=N(m,\delta)\geq1$ and $N^{*}=N^{*}(m,\delta)$
such that 
if

\begin{equation}
                                              \label{12.23.4}
\nu\geq N^{*} \|\eta'\|+N (\|\eta''\|+\|\eta'\|^{2}) 
\end{equation}
and $N^{*}h\leq 1$, then  on $Q^{o}$
for any $\alpha=(a,b,c)\in A$ we have

$$
2L_{\nu}[a_{k}\Delta_{k}+b_{k}\delta_{k}]u\geq
-[a_{k}\Delta_{k}+b_{k}\delta_{k}]V_{\nu}
$$

\begin{equation}
                                            \label{9.20.3}
-
(N\nu^{2}+N^{*}\nu)\bar{W}_{r}+(\nu/2)\zeta 
a_{k}u_{kr}^{2}.
\end{equation} 
Furthermore, $N^{*}=0$ if $b\equiv0$.

\end{lemma}

Proof. Since, for each $k$,
 the operators $a_{k}\Delta_{k}+b_{k}\delta_{k}$
respect the maximum principle,
it follows by  Lemma 4.2 of \cite{Kr2} that

$$
u_{rr}^{-}(a_{k}\Delta_{k}+b_{k}\delta_{k})u_{rr} 
\geq-u_{rr}^{-}[a_{k}\Delta_{k}+b_{k}\delta_{k}]u_{rr}^{-}.
$$
Hence,
$$
I:=\zeta^{2}u^{-}_{rr}\Delta_{r} 
[a_{k}\Delta_{k}+b_{k}\delta_{k}]u=
\zeta^{2}u^{-}_{rr}  
[a_{k}\Delta_{k}+b_{k}\delta_{k}]u_{rr}
$$

$$
\geq-\zeta^{2}u^{-}_{rr}  
[a_{k}\Delta_{k}+b_{k}\delta_{k}]u^{-}_{rr},
$$
which by Lemma \ref{lemma 9.19.1} and the observation that

$$
\zeta (u_{rr}^{-})^{2} \leq N \zeta a_{k}u_{kr}^{2},
$$
for $N^{*}h\leq 1$ yields

$$
2I\geq -[a_{k}\Delta_{k}+b_{k}\delta_{k}]
(\zeta^{2}(u^{-}_{rr})^{2})-N(N^{*}|\eta'|+|\eta''|
+|\eta'|^{2})\zeta a_{k}u_{kr}^{2}
$$

\begin{equation}
                                             \label{12.22.2}
-(N^{*}|\eta'|^{2}+N|\eta'|^{4})\bar{W}_{r}.
\end{equation}

Furthermore, by Lemma \ref{lemma 9.19.1}
$$
-2\nu\zeta u_{r}\delta_{r}[a_{k}\Delta_{k}+b_{k}\delta_{k}]
u=-2\nu\zeta u_{r} [a_{k}\Delta_{k}+b_{k}\delta_{k}]
u_{r}
$$

$$
\geq -[a_{k}\Delta_{k}+b_{k}\delta_{k}](\nu\zeta u^{2}_{r})
+\nu\zeta a_{k} u^{2}_{kr}-N\nu(|\eta''|+|\eta'|^{2}
+N^{*}|\eta'|)\bar{W}_{r}.
$$
By combining this with \eqref{12.22.2} and recalling
\eqref{9.21.1} we find

 $$
2L_{\nu}[a_{k}\Delta_{k}+b_{k}\delta_{k}]u
\geq -[a_{k}\Delta_{k}+b_{k}\delta_{k}]V_{\nu}
$$

$$
+( \nu-N^{*}_{1}|\eta'|-N_{1}|\eta''|-N_{1}|\eta'|^{2})\zeta 
a_{k}u_{kr}^{2}
$$

\begin{equation}
                                            \label{9.20.2}
-[N ^{*}(|\eta'|^{2}+\nu|\eta'|)+N(|\eta'|^{4}
+  \nu |\eta''|+\nu |\eta'|^{2})
 ] \bar{W}_{r}.
\end{equation}  
We may assume that $N_{1}\geq1$ and then,
if

$$
\nu\geq2 N^{*}_{1}|\eta'|+2N_{1}(|\eta''|+|\eta'|^{2}),
$$
we have that $|\eta'|^{2}\leq\nu$, $|\eta''|\leq\nu$, and

$$
|\eta'|^{4}
+  \nu |\eta''|+\nu |\eta'|^{2}\leq 3\nu^{2}.
$$
Also

$$
N ^{*}(|\eta'|^{2}+\nu|\eta'|)\leq
N ^{*}(\nu+\nu^{3/2} )\leq N ^{*}(\nu+\nu^{2} )
\leq N^{*}\nu+N\nu^{2}.
$$
After that \eqref{9.20.2} clearly yields \eqref{9.20.3}
and the lemma is proved.

{\bf Proof of Theorem \ref{theorem 9.18.1}}.
Denote by $N_{0}$ and $N^{*}_{0}$ 
the constants $N$ and $N^{*}$ in Lemma \ref{lemma 9.20.1}
and take and fix a  $\nu$ satisfying
\eqref{12.23.4} (with $N_{0}$ and $N^{*}_{0}$ in place
of $N$ and $N^{*}$).

  Observe that
\eqref{9.18.3} obviously holds on $Q\setminus Q^{o}$
and we only need to prove it on $Q^{o}$. Also notice that
$$
|u_{rr}|=h^{-1}|u_{r}+u_{-r}|\leq 2h^{-1}\bar{W}_{r}^{1/2},
$$
which shows that \eqref{9.18.3} holds if $h\geq\nu^{-1/2}$
or if $N^{*}_{0}h\geq1$.
Therefore below we
  assume that
\begin{equation}
                                      \label{12.23.1}
h\leq\nu^{-1/2},\quad N^{*}_{0}h\leq1.
\end{equation}  

Introduce
$$
D^{o}=\{x\in Q^{o}:\zeta(x)u_{rr}^{-}(x)\geq\nu h u_{ r}(x)\}.
$$
If $x\in  Q^{o}\setminus D^{o}$, then 
$\zeta(x)u_{rr}^{-}\leq\nu h u_{ r}(x)$ and
 \eqref{9.18.3}  holds at $x$ in light of 
\eqref{12.23.1}. 
Therefore, we in the remaining part of the proof
we   concentrate on establishing  \eqref{9.18.3} 
for $x\in D^{o}$,
assuming, of course, that  $D^{o}\ne\emptyset$.

Denote 

$$
D=D^{o}\cup\{x+h\Lambda:x\in D^{o}\}\quad (\subset Q).
$$
If $V_{\nu}$ in $D^{o}$ is less than its maximum
over $D\setminus D^{o}$, then in $D^{o}$
$$
\zeta^{2}[u_{rr}^{-}]^{2}
\leq V_{\nu}\leq \max[
\max_{Q\setminus Q^{o}}\zeta^{2}[u_{rr}^{-}]^{2},
\max_{ Q^{o}\setminus D^{o}}\zeta^{2}[u_{rr}^{-}]^{2}]
+\nu \bar{W}_{r},
$$
where the maximums are less than the right-hand side of
\eqref{9.18.3} by the above. Hence, in the rest of the proof
 we  consider the case that
the maximum over $D$ of $V_{\nu}$ is attained at a point
$x_{0}\in D^{o}$.

Notice
that if a function $\phi(x)$ is such that
 $\phi(x)\leq\phi(x_{0})$
for $x\in x_{0}+h\Lambda$, then at $x_{0}$

$$
h^{2}L_{\nu}\phi(x_{0})=\zeta [
\phi(x_{0}+h\ell_{r})(\zeta u_{rr}^{-} -\nu h u_{r})
+  \phi(x_{0}-h\ell_{r})\zeta u_{rr}^{-}]
$$

$$
-\zeta[2\zeta u_{rr}^{-}-\nu h u_{r} ]\phi
\leq
\zeta [
(\zeta u_{rr}^{-} -\nu h u_{r})\phi 
+ \zeta u_{rr}^{-} \phi ]-\zeta[2\zeta u_{rr}^{-}-\nu h u_{r} ]\phi,
$$
where the last expression is zero. Thus
$$
L_{\nu}\phi(x_{0})\leq0,
$$
which in the terminology from \cite{Kr2} means that
$L_{\nu}$ respects the maximum principle.

Next, we can find an  
$\bar{\alpha} =(\bar{a},\bar{b},\bar{c})\in A$
  such that
$$
\bar{a}_{k}\Delta_{ k}u(x_{0})
+\bar{b}_{k}\delta_{ k}u(x_{0})
-\bar{c}u(x_{0})+P(\bar{\alpha},x_{0})=
\cP[u](x_{0})=0.
$$
Since $\cP[u]\leq0$ in $Q $, we have that 
$$
\phi(x):=
\bar{a}_{k}\Delta_{ k}u(x )
+\bar{b}_{k}\delta_{ k}u(x )
-\bar{c} u(x )+f(\bar{\alpha},x )\leq0
$$
for $x\in x_{0}+h\Lambda$. Hence,  
$$
0\geq  2L_{\nu}\phi(x_{0})=
\bar{a}_{k}2L_{\nu}\Delta_{ k}u(x_{0} )
+\bar{b}_{k}2L_{\nu}\delta_{ k}u(x_{0})
$$

$$
- \bar{c} 2L_{\nu}u(x_{0})+
2L_{\nu}f(\bar{\alpha},\cdot )(x_{0} ) ,
$$
which owing to \eqref{9.21.1} and \eqref{9.20.3}
yields
$$
0\leq[\bar{a}_{k}\Delta_{k}+\bar{b}_{k}\delta_{k}-
2\bar{c} 
]V_{\nu}(x_{0})- (\nu/2)\zeta 
a_{k}u_{kr}^{2}(x_{0})
$$

$$
+(N\nu^{2}+N^{*}\nu)\bar{W}_{r}
- 2L_{\nu}f(\bar{\alpha},\cdot )(x_{0} ) .
$$
Here the last term is dominated by
$$
 K_{2}\zeta^{2}u_{rr}^{-}(x_{0})+\nu |u_{r}(x_{0})| K_{1} 
$$

$$
\leq  N\nu^{-1}K_{2}^{2}+(\nu/4)\zeta a_{k}u_{kr}^{2}
(x_{0})
+ K_{1}^{2}+\nu^{2}\bar{W}_{r} .
$$
Furthermore, by the maximum  principle
$$
[\bar{a}_{k}\Delta_{k}+\bar{b}_{k}\delta_{k}-
2\bar{c} 
]V_{\nu}(x_{0})\leq0,
$$
since $V_{\nu}\geq0$ attains its maximum at $x_{0}$.

We now conclude that
$$
(\nu/4)\zeta a_{k}u_{kr}^{2}
(x_{0})\leq (N\nu^{2}+N^{*}\nu)\bar{W}_{r}
+N\nu^{-1}K_{2}^{2} 
+ K_{1}^{2}  ,
$$
which implies that in $D^{o}$
$$
\zeta^{2}(u_{rr}^{-})^{2}\leq
V_{\nu}(x_{0})\leq N\zeta a_{k}u_{kr}^{2}
(x_{0})+\nu\bar{W}_{r}
$$

$$
\leq (N\nu+N^{*})\bar{W}_{r}+N\nu^{-2}K_{2}^{2} 
+ \nu^{-1}K_{1}^{2}.
$$
Thus, estimate \eqref{9.18.3} holds on $D^{o}$
as well and this proves the theorem.

\mysection{A model cut-off equation}  

                                        \label{section 10.18.2}

We will work in the setting of
 Section \ref{section 9.22.1}. However now  $h>0$ is not fixed.
Take a function $\cH(u,x)$, where $x\in\bR^{d}$,
$u=(u',u'')\in\bR^{1+m'+2m}$. 

\begin{assumption}
                                       \label{assumption 9.23.01}
(i) The function $\cH$ is Lipschitz continuous in $u$
for every $x$, and at all points of differentiability
of $\cH$ with respect to $u$   we have
$$
\delta\leq \cH_{u''_{k}}\leq\delta^{-1},\quad
  k=\pm1,...,\pm m,
\quad \delta\leq- \cH_{u' _{0}}\leq\delta^{-1},
$$
$$
|\cH_{u'_{k}}|\leq\delta^{-1},\quad k= 1,..., m';
$$

(ii) The number $\bar{\cH}=\sup_{x}|\cH(0,0,x)|$ is finite;

(iii) The function $\cH$ is locally Lipschitz continuous in $x$
for every $u$ and there exists a constant $N'$
such that at all points of differentiability
of $\cH$ with respect to $x$   we have
$$
|\cH_{x_{i}}(u,x)|\leq N'(1+|u|),\quad\forall i;
$$

(iv) We have ${\rm Span}\,(l_{1},...,l_{m})=\bR^{d}$.

\end{assumption}

Define 
$$
\cP (u',u'',x)=\cP (u' ,u'')=
2\delta ^{-1}\sum_{k}(u''_{k})^{+}
-(\delta /2)\sum_{k}(u''_{k})^{-} 
$$
$$
+2\delta ^{-1}\sum_{k\geq1}|u'_{k}|
-(\delta/2) (u'_{0})^{+}+2\delta ^{-1}
(u'_{0})^{-}
$$
 
\begin{equation}
                                               \label{9.22.2}
=\max_{\substack{\delta /2\leq a_{k}\leq
2/\delta \\|k|=1,...,m} }
\max_{\substack{ |b_{k}|\leq
2/\delta \\|k|=1,...,m} }
\max_{\delta/2  \leq c\leq
2/\delta }\big[\sum_{|i|=1}^{m} a_{i}
u''_{i}+\sum_{ i =1}^{m'}b_{i}u'_{i} -cu'_{0}\big].
\end{equation}

For  functions $v(x)$   introduce
$$
H[v](x)=\cH(v(x),\partial v(x),\partial^{2}v(x),x)
$$
whenever and wherever it makes sense, where
$$
\partial v=(  
v_{(\ell_{1})},...,v_{(\ell_{ m'})}),
$$
$$
\partial^{2} v=( v_{(\ell_{-m})(\ell_{-m})},...,v_{(\ell_{-1})(\ell_{-1})},
v_{(\ell_{1})(\ell_{1})},...,v_{(\ell_{ m})(\ell_{ m})}),
$$
and $v_{(\ell)}=\ell_{i}v_{x_{i}}$, $v_{(\ell)
(\ell)}=\ell_{i}\ell_{j}v_{x_{i}x_{j}}$. Similarly,
$$
P[u](x)=\cP(u(x),\partial u(x),\partial^{2}u(x) ).
$$

Let $\Omega$ be a bounded $C^{2}$ domain in $\bR^{d}$,
$g\in C^{1,1}(\Omega)$,
and let $K\geq0$ be a finite number.
\begin{theorem}
                                       \label{theorem 9.14.1}
In addition to
Assumption \ref{assumption 9.23.01} suppose
  that $\pm e_{i},\pm (e_{i}+ e_{j}),
  e_{i} -e_{j}\in \Lambda$, $i,j=1,..,d$,
 were  $e_{1},...,e_{d}$
is the standard orthonormal basis in $\bR^{d}$  
and  assume that all vectors in $\Lambda$ have rational 
coordinates.
Then there exists
a unique $v\in C^{0,1}(\bar{\Omega})\cap C^{1,1}_{loc}(\Omega)$
such that $v=g$ on $\partial \Omega$ and
\begin{equation}
                                              \label{9.14.2}
H _{K}[v]=0
\end{equation}
(a.e.) in $\Omega$, where
$$
 H _{K}[v]=
\max( H[v],P [v]-K).
$$
  Furthermore,  
\begin{equation}
                                             \label{1.14.3}
|v| ,  |D_{i}v| ,
  \rho|D_{ij} v|\leq N(\bar{\cH}+K+\|g\|_{C^{1,1}(\Omega)})
\end{equation}
 in $\Omega$ (a.e.)
for all $i,j$,
where 
 $N$ is a   constant depending only
on $\Omega$, $\{\ell_{1},...,\ell_{m}\}$, $d$, and $\delta$
(but not on $N'$).

\end{theorem}

To prove the theorem,
we are going to use finite-difference approximations
of the operators $H[v]$ and $P[v]$.
For $h>0$ introduce
$$
P_{h}[v](x)=\cP(v(x),\delta_{h}v(x),\delta_{h}^{2}v(x)),
$$
where naturally  
$$
\delta_{h} u=( \delta_{h,1}u ,...,
\delta_{h,m'}u ),
$$
$$
\delta^{2}_{h} u=(\Delta_{h,-m}u ,...,
\Delta_{h,-1}u ,\Delta_{h,1}u ,...,
\Delta_{h,m}u).
$$
Similarly we introduce $H_{h}$ and $H_{K,h}$.

\begin{lemma}
                                          \label{lemma 9.14.1}
Under Assumptions \ref{assumption 9.23.01} (i), (ii)
\begin{equation}
                                            \label{9.14.1}
\cH \leq  \cP -(\delta /2)
 \big(\sum_{k}|u''_{k}| +\sum_{k}|u'_{k}|\big)+
\bar{\cH}.
\end{equation}
\end{lemma}

Proof. Basically,
 we just repeat part of the proof of Lemma \ref{lemma 9.29.2}.
From Hadamard's formula
$$
H (u' ,u'',x)-H (0,0,x)
=u''_{k}\int_{0}^{1}H _{u''_{k}}
(tu' ,tu'',x)\,dt
$$
$$
+\sum_{k\geq1}
u' _{k}\int_{0}^{1}H _{u' _{k}}
(tu' ,tu'',x)\,dt+
u' _{0}\int_{0}^{1}H _{u' _{0}}
(tu' ,tu'',x)\,dt
$$
we obtain
$$
H (u' ,u'',x)-H (0,0,x)\leq
\delta ^{-1}\sum_{k}(u''_{k})_{+}-
\delta  \sum_{k}(u''_{k})_{-}
$$
$$
+\delta ^{-1}\sum_{k\geq1}|u' _{k}|
-\delta (u'_{0})_{+}+\delta ^{-1}
(u'_{0})_{-}
$$
$$
=P (u' ,u'')-\delta ^{-1}\sum_{k} (u''_{k})_{+}
-(\delta /2)\sum_{k}(u''_{k})_{-} 
-\delta ^{-1}\sum_{k\geq1}|u'_{k}|
$$
$$
-(\delta/2) (u'_{0})_{+}+\delta ^{-1}
(u'_{0})_{-}
$$
and \eqref{9.14.1} follows. The lemma is proved.

Introduce $B$ as the smallest closed
 ball containing $\Lambda$ and set
$$
\Omega_{h}=\{x\in \Omega:x+hB\subset \Omega\}
=\{x:\rho(x)\geq \lambda h\},
$$
where $\lambda$ is the radius of $B$.

Consider the equation
\begin{equation}
                                              \label{2.25.3}
H _{K,h}[v]=0\quad
\text{in}\quad \Omega_{h}   
\end{equation}
 with boundary condition
\begin{equation}
                                              \label{2.25.4}
v=g\quad\text{on}\quad \Omega\setminus \Omega_{h} .
\end{equation}

It is a rather simple fact that for sufficiently small
$h>0$ there exists
a unique bounded solution $v=v_{ h}$ of 
\eqref{2.25.3}--\eqref{2.25.4} (see, for instance,
\cite{KT} or Theorem 8.2 in \cite{Kr11} or else
Theorem 2.2 in \cite{Kr12.1}). By the way,
we do not include $K$ in the notation $v_{ h}$
since $K$ is a fixed number.

Below by $h_{0}$ and $N$ with occasional indices
we denote various (finite) constants depending only
on $\Omega$, $\{l_{1},...,l_{m}\}$, $d$, and $\delta$.

In the following   lemma  the additional assumption of
Theorem \ref{theorem 9.14.1} concerning the $e_{i}$'s
and the $\ell_{k}$'s
is not used.
\begin{lemma}
                                         \label{lemma 9.14.2}
Under Assumptions \ref{assumption 9.23.01} (i), (ii), (iv)
there are constants $h_{0}>0$ and
$N$ 
 such that for all $h\in(0,h_{0}]$
and $|r|=1,...,m$
\begin{equation}
                                              \label{1.14.1}
|v_{ h}-g|\leq N(\bar{\cH}+K+\|g\|_{C^{1,1}(\Omega)})\rho,
\end{equation}
\begin{equation}
                                              \label{1.14.2}
  |\delta_{h,r}v_{ h}|
\leq N(\bar{\cH}+K+\|g\|_{C^{1,1}(\Omega)})
\end{equation}
on $\Omega$.
\end{lemma}

Proof. Introduce
$$
\cH_{K}=\max(\cH,\cP-K).
$$
Clearly, $\cH_{K}$ satisfies Assumption \ref{assumption 9.23.01}
 with $\delta/2$ in place of $\delta$. Therefore, by Hadamard's
formula there exist functions $a_{k},b_{k}$, $k=\pm1,...,\pm m$,
and $c$ such that
\begin{equation}
                                                \label{10.18.4}
\delta/2\leq a_{k}\leq2\delta^{-1},\quad|b_{k}|\leq2\delta^{-1},
\quad \delta/2\leq c\leq2\delta^{-1}
\end{equation}
and in $\Omega_{h}$ we have
$$
-\cH_{K}[0]=H_{K,h}[v_{h}]-\cH_{K}[0]=
a_{k}\Delta_{h,k}v_{h}+b_{k}\delta_{h,k}v_{h}-
cv_{h}
$$
$$
=
a_{k}\Delta_{h,k}(v_{h}-g)+b_{k}\delta_{h,k}(v_{h}-g)-
c(v_{h}-g)+f,
$$
where $$
f=
a_{k}\Delta_{h,k}g+b_{k}\delta_{h,k}g-
cg.
$$
After that   \eqref{1.14.1}
is proved by using the barrier function $\Phi$ from Lemma
 2.4 of \cite{Kr12.1} (cf.  Lemma 2.5 in \cite{Kr12.1}).
 It implies that
\begin{equation}
                                           \label{9.16.3}
|v_{ h}-g|\leq N_{1}(\bar{\cH}+K
+\|g\|_{C^{1,1}(\Omega)})h \quad \text{on}\quad
\Omega\setminus\Omega_{3h}
\end{equation}
with a constant $N_{1}$
 independent of $h$.

To prove   \eqref{1.14.2},
fix an $r$ and define
$$
Q^{o} =\{x\in\Omega_{2h}:(\delta/2)|\delta_{h,r}v_{h}|
\geq \bar{\cH}+K\}.
$$
If $Q^{o}=\emptyset$, then $(\delta/2)|\delta_{h,r}v_{h}|
\leq \bar{\cH}+K$ in $\Omega_{2h}$, and by virtue of 
\eqref{9.16.3}, 
$$
|\delta_{h,r}(v_{h}-g)| \leq 2N_{1}(\bar{\cH}+K
+\|g\|_{C^{1,1}(\Omega)})
$$
in $\Omega\setminus\Omega_{2h}$. In that case
  \eqref{1.14.2} obviously holds.

Therefore, we assume that $Q^{o} \ne\emptyset$ and owing to
Lemma \ref{lemma 9.14.1} conclude that 
\begin{equation}
                                              \label{9.16.7}
P_{h} [v_{ h}]=K 
\end{equation}
in $Q^{o}$. Furthermore, \eqref{2.25.3} implies that 
\begin{equation}
                                              \label{9.22.3}
P_{h}[v_{h}]\leq K
\end{equation}
in $\Omega_{h}$.  

Now use again the mean value theorem to conclude that
$$
\delta_{h,r}P_{h}[v_{h}]=
a_{k}\Delta_{h,k}(\delta_{h,r}v_{h})
+b_{k}\delta_{h,k}(\delta_{h,r}v_{h})-c(\delta_{h,r}v_{h})
$$
for some functions $a_{k}(x),b_{k}(x),c(x)$ satisfying 
\eqref{10.18.4}.
Furthermore, $\delta_{h,r}P_{h}[v_{h}]\leq 0$ in $Q^{o}$
owing to \eqref{9.16.7} and \eqref{9.22.3}, that is  in $Q^{o}$
$$
a_{k}\Delta_{h,k}(\delta_{h,r}v_{h})
+b_{k}\delta_{h,k}(\delta_{h,r}v_{h})-c(\delta_{h,r}v_{h})
\leq0.
$$

For small enough $h_{0}$  the operator
$a_{k}\Delta_{h,k} 
+b_{k}\delta_{h,k} -c$ with $h\in(0,h_{0}]$
respects the maximum principle
and therefore in $Q^{o}$ (see Theorem 2.2 in \cite{Kr12.1})
\begin{equation}
                                              \label{9.22.4}
 (\delta_{h,r}v_{h})_{+}\leq\sup_{\Omega\setminus Q^{o}}
(\delta_{h,r}v_{h})_{+}.
\end{equation}
Notice that if $x\in \Omega\setminus Q^{o}$,
then either $x\not\in\Omega_{2h}$, in which case
\eqref{1.14.2} holds by the above,  or else $x\in
\Omega_{2h}$ but $(\delta/2)|\delta_{h,r}v_{h}|
\leq \bar{\cH}+K$. It follows that the left-hand side
of \eqref{9.22.4} is dominated by the right-hand side of
\eqref{1.14.2},
  if $h\in(0,h_{0}]$
and $h_{0}>0$  is sufficiently small.

Thus, in all cases   
$$
(\delta_{h,r}v_{h})_{+}\leq N(\bar{\cH}+K
+\|g\|_{C^{1,1}(\Omega)})
$$
on $\Omega$. Upon replacing here $r$ with $-r$, we get
$$
T_{h,-\ell_{r}}(\delta_{h,r}v_{h})_{-}\leq N(\bar{\cH}+K
+\|g\|_{C^{1,1}(\Omega)}) 
$$ 
in $\Omega_{h}$, which after being combined with the previous
estimate proves \eqref{1.14.2} in $\Omega_{h}$.
In $\Omega\setminus\Omega_{h}$ estimate
\eqref{1.14.2} has been established above. The lemma is
proved.

\begin{lemma}
                                          \label{lemma 9.16.1}
Suppose that
Assumptions \ref{assumption 9.23.01} (i), (ii), (iv) are satisfied.
Assume also that all vectors in $\Lambda$ have rational 
coordinates. Then there are constants $h_{0}>0$ and
$N$ 
such that for all $h\in(0,h_{0}]$
and $|r|=1,...,m$
\begin{equation}
                                              \label{9.16.6}
 (\rho-6\lambda h) |\Delta_{h,r}  v_{ h}|\leq N(\bar{\cH}+K
+\|g\|_{C^{1,1}(\Omega)})
\end{equation}
on $\bR^{d}$ (we remind the reader that
$\lambda$ is the radius of $B$).
\end{lemma}

Proof. Clearly, the assertion of the lemma
would follow if we can prove that \eqref{9.16.6}
holds
on $y+h\Lambda_{ \infty}$ for any $y\in\bR^{d}$ with a constant
$N$ independent of $h$  and $y$. Without losing generality
we concentrate on $y=0$. Then for a fixed $r$ define
$$
Q^{o}:=\{x\in (h\Lambda_{ \infty})\cap\Omega_{3h}
 :(\delta/2)  |\Delta_{h,r}
v_{ h}(x)|\geq  \bar{\cH}+ K \}.
$$
If $x\in h\Lambda_{\infty}$
is such that $x\not\in Q^{o}$, then either
$x\not\in \Omega_{3h}$, so that $\rho(x)\leq 
3\lambda h$ and \eqref{9.16.6} holds, or else
$x \in \Omega_{3h}$ but $(\delta/2)  |\Delta_{h,r}
v_{ h}(x)|\leq  \bar{\cH}+ K$, in which case
\eqref{9.16.6} holds again.

Thus we need only prove \eqref{9.16.6} on $Q^{o}$
assuming, of course, that $Q^{o}
\ne\emptyset$.
Then  define
$$
Q=\{x+h\Lambda :x\in Q^{o}\}   .
$$
Observe that
$Q$ is a finite set since $\ell_{k}$ have rational coordinates and 
 there is a number
$M$ such that the coordinates of all points in 
$M\Lambda_{1,\infty}$
are integers and
the number of points with integral coordinates
lying in a  bounded domain
is finite.
 
Next by Lemma \ref{lemma 9.14.1} we have
that \eqref{9.16.7} holds in $Q^{o}$ and
\eqref{9.22.3} holds in $Q\setminus Q^{o}$.

To proceed further observe 
 a standard fact that there are constants $\mu_{0}>0$
and $N\in[0,\infty)$
depending only on $\Omega$  
such that
 for any $\mu\in(0,\mu_{0}]$ there exists an $\eta_{\mu}
\in C^{\infty}_{0}(\Omega)$ satisfying
$$
\eta_{\mu}=1\quad\text{on}\quad \Omega_{2\mu},\quad
\eta_{\mu}=0\quad\text{outside}\quad \Omega_{\mu},
$$
\begin{equation}
                                                    \label{9.22.6}
|\eta_{\mu}|\leq1,\quad
|D\eta_{\mu}|\leq N/\mu,\quad|D^{2}\eta_{\mu}|\leq N/\mu^{2}.
\end{equation}
By Theorem \ref{theorem 9.18.1} and Lemma \ref{lemma 9.14.2}
 there are  constants $N$ and $h_{0}>0$
such that, for any number
$\nu$ satisfying
$$
\nu\geq N(\|\eta'_{\mu}\|_{h}+\|\eta'_{\mu}\|_{h}^{2}+
\|\eta''_{\mu}\|_{h}),  
$$
we have in $Q^{o}$ that
$$
\eta_{\mu}^{4}[ ( \Delta_{r}v_{h})^{-}]^{2}\leq
\max_{Q\setminus Q^{o}}\eta_{\mu}^{4}[ ( \Delta_{r}v_{h})^{-}]^{2}
+ N(\nu+1) (\bar{\cH}+ K+\|g\|_{C^{1,1}(\Omega)})^{2}
$$
 if $h\in(0,h_{0}]$. In light of \eqref{9.22.6}
one can take
$\nu=N\mu^{-2}$ for an appropriate $N$ and then
$$
 \eta_{\mu}^{4}[ ( \Delta_{r}v_{h}(x))^{-}]^{2}\leq
\max_{y\in Q\setminus Q^{o}}\eta_{\mu}^{4}[ ( \Delta_{r}v_{h}
(y))^{-}]^{2}
+ N \mu ^{-2}(\bar{\cH}+ K+\|g\|_{C^{1,1}(\Omega)})^{2}
$$
for $x\in Q^{o}$. We will only 
concentrate on $\mu\geq 3h$, when $\eta_{\mu}=0$
outside $\Omega_{3h}$. In that case, for any $y\in Q\setminus
Q^{o}$, either $y\not\in\Omega_{3h}$ implying that
$$
\eta_{\mu}^{4}[ ( \Delta_{r}v_{h})^{-}]^{2}(y)=0
$$ 
or else $y \in\Omega_{3h}\cap(h\Lambda_{\infty})$ but
$$
(\delta/2)  |\Delta_{h,r}
v_{ h}(y)|\leq  \bar{\cH}+ K.
$$

It follows that as long as $h\in(0,h_{0}]$,
$x\in Q^{o}$,
and  $\mu\geq 3h$ we have
\begin{equation}
                                                   \label{9.22.7}
\eta_{\mu}^{4}[ ( \Delta_{r}v_{h})^{-}(x)]^{2}\leq
 N \mu ^{-2}(\bar{\cH}+ K+\|g\|_{C^{1,1}(\Omega)})^{2}.
\end{equation}

If $x$ is such that $\rho(x)\geq 6\lambda h$, take
$\mu=\mu_{0}\wedge(\rho(x)/(2\lambda))$, which is bigger than
$3h$ provided that $h\leq\mu_{0}/3$. In that case also
$$
\rho(x)=2\lambda [\rho(x)/(2\lambda)]\geq 2\lambda \mu,
$$ 
so that $\eta_{\mu}(x)=1$ and we conclude from
\eqref{9.22.7} that
$$
\rho(x)( \Delta_{r}v_{h})^{-}(x)
\leq N(\bar{\cH}+ K+\|g\|_{C^{1,1}(\Omega)}), 
$$
\begin{equation}
                                                   \label{9.22.8}
(\rho(x)-6\lambda h)( \Delta_{r}v_{h})^{-}(x)
\leq N(\bar{\cH}+ K+\|g\|_{C^{1,1}(\Omega)})
\end{equation}
for $x\in Q^{o}$ such that $\rho(x)\geq 6\lambda h$.
However, the second relation in 
\eqref{9.22.8} is obvious for $\rho(x)\leq 6\lambda h$.
 
As a result of all the above arguments we see that
\eqref{9.22.8} holds in $h\Lambda_{\infty}$ for any $r$
whenever $h\in(0,h_{0}]$.

Finally, since $P_{h}[v_{h}]\leq K$ in $\Omega_{h}$ we have
that
$$
2\delta ^{-1}\sum_{r}(\Delta_{r}v_{h})_{+}
\leq(\delta /2)\sum_{r}(\Delta_{r}v_{h})_{-} 
$$
$$
-2\delta ^{-1}\sum_{r\geq1}|\delta_{r}v_{h}|
+(\delta/2) (v_{h})_{+}-2\delta ^{-1}
(v_{h})_{-}+K,
$$
which after being multiplied by $\rho-6h$ along with
\eqref{9.22.8} and Lemma \ref{lemma 9.14.2} leads
to \eqref{9.16.6} on $h\Lambda_{\infty}$.
As  is explained in the beginning of the proof,
this finishes proving the lemma.

{\bf Proof of Theorem \ref{theorem 9.14.1}}.
Owing to Assumption \ref{assumption 9.23.01} (iii),
by Corollary 2.7 of \cite{Kr12.1}, which is applicable
in light of Lemmas \ref{lemma 9.14.2} and \ref{lemma 9.16.1},
there exists a constant $M$ such that for all sufficiently
small $h$ and $x,y\in\bR^{d}$ we have
$$
|v_{h}(x)-v_{h}(y)|\leq M(|x-y|+h).
$$
Here Assumptions \ref{assumption 9.23.01}(iii) plays a  
crucial role.

After that our  theorem is proved in exactly the same way
as
Theorem 8.7 of \cite{Kr11} on the basis
of Lemmas \ref{lemma 9.14.2} 
and \ref{lemma 9.16.1} 
and the fact that the derivatives of $v$ are weak limits
of finite differences of $v_{h}$ as $h\downarrow0$
(see the proof of Theorem 8.7 of \cite{Kr11}).
One also uses the fact 
that there are sufficiently many pure
 second order derivatives in directions
of the $l_{i}$'s  to conclude from their boundedness
that the Hessian of $v$ is bounded.
The theorem is proved.

\mysection{Proof of Theorem \protect\ref{theorem 10.5.1}}
                                      \label{section 12.13.4}

The functions $\cH$ from Section \ref{section 10.18.1} and $\cP$ 
from Section \ref{section 2.5.1}
are instances of $\cH$ and $\cP$ from Section
\ref{section 10.18.2}. To see this, of course, one
has to change the constant $\delta$ in 
Section \ref{section 10.18.2}  and renumber the $l_{i}$'s 
in Section \ref{section 2.5.1}. 
We also take into account that $\hat{\delta}\leq\delta/4$
which allows us to match \eqref{5.10.1} and \eqref{5.10.2}
with the requirements of Assumption \ref{assumption 9.23.01} (i).
Furthermore,  
$\bar{\cH}=\bar{H}$.
Therefore, Theorem \ref{theorem 9.14.1} is applicable and
  yields a unique $v\in C^{0,1}(\bar{\Omega})
\cap C^{1,1}_{loc}(\Omega)$ such that $v=g$ on $\partial
\Omega$, estimates \eqref{1.14.3}, that is \eqref{1.13.1},  
 hold true, and
$$
\max[\cH(v,Dv,v_{(l_{1})(l_{1})},...,v_{(l_{m})(l_{m})},x),
\cP(v,Dv,v_{(l_{1})(l_{1})},...,v_{(l_{m})(l_{m})})-K]=0
$$
in $\Omega$ (a.s.). In light of the construction of
 $\cH$ and $\cF$ in Section \ref{section 10.18.1} this equation
coincides with \eqref{9.23.2}, so that
  the only remaining assertions of Theorem 
\ref{theorem 10.5.1} to prove are that for $p>d$
\begin{equation}
                                                \label{10.19.1}
\|v\|_{W^{2}_{p}(\Omega)}\leq
N_{p}(\bar{H}+K+\|g\|_{W^{2}_{p}(\Omega)})
\end{equation}
and estimate \eqref{2.28.1} holds. The latter
follows from other assertions
of Theorem \ref{theorem 10.5.1} by Remark
\ref{remark 2.28.1}, so that we may concentrate on 
\eqref{10.19.1}.

Observe that
$$
\max(H(u,x),P(u)-K)=P(u)+G(u,x),
$$
where $G(u,x)= (H(u,x)-P(u)+K)_{+}-K$ and, owing to
condition \eqref{3.6.4},
$G(u,x)=-K$ if
$$
\kappa\big(\sum_{i,j}|u_{ij}|+\sum_{i}|u_{i}|\big)
\geq\bar{H}+K.
$$
If the opposite inequality holds, then
\begin{equation}
                                            \label{1.14.4}
|G(u,x)|\leq | H(u,x)-H(0,x)|+|P(u)|+\bar{H}+K
\leq N(\bar{H}+K),
\end{equation}
where $N$ depends only on $\delta$ and $d$.
It follows that the inequality between
the extreme terms in \eqref{1.14.4} holds for all $u$ and $x$.
This allows us to apply Theorem 1.2 of \cite{DKL}
and shows that \eqref{10.19.1} holds if
$v\in W^{2}_{p}(\Omega)$ or if $w:=v-g\in W^{2}_{p}(\Omega)$.
In light of \eqref{1.13.1} it suffices to show that
$w \in W^{2}_{p}(D\cap\Omega)$, where $D$ is a neighborhood of
$\partial\Omega$.

To prove the latter  we use the fact that, according to
\eqref{10.5.2}, in a neighborhood $D$ of
$\partial\Omega$ intersected with $\Omega$ the function $w$
satisfies the equation
\begin{equation}
                                          \label{10.19.2}
\max(\Delta u-u+\Delta g-g,P[u+g]-K)=0,
\end{equation}
in which the left-hand side is given by a convex
function of $u$ and its derivatives.
We may certainly assume that $D\in C^{2}$ and then, by relying
on \eqref{1.13.1},
find a $\zeta\in C^{\infty}_{0}(\Omega)$
such that  $\zeta w=w$ on $\partial (D\cap\Omega)$
and $\zeta w\in W^{2}_{p}(D\cap\Omega)$. Then
due to Theorem 1.2 of \cite{DKL} equation \eqref{10.19.2}
with boundary condition $u-\zeta w
\in \WO^{2}_{p}(D\cap\Omega)$ 
has a  (unique) solution $u\in W^{2}_{p}
(D\cap\Omega)$.

By uniqueness of $W^{2}_{d,loc}
(D\cap\Omega)\cap C(\overline{D\cap\Omega})$-solutions
we obtain $w=u\in W^{2}_{p}
(D\cap\Omega)$ and the theorem is proved.

\end{document}